\documentclass[pdflatex,sn-mathphys-num]{sn-jnl}

\usepackage{graphicx}%
\usepackage{multirow}%
\usepackage{amsmath,amssymb,amsfonts}%
\usepackage{amsthm,bm}%
\usepackage{mathrsfs}%
\usepackage[title]{appendix}%
\usepackage{xcolor}%
\usepackage{textcomp}%
\usepackage{manyfoot}%
\usepackage{booktabs}%
\usepackage{algorithm}%
\usepackage{algorithmicx}%
\usepackage{algpseudocode}%
\usepackage{listings}%
\usepackage{mathabx}
\usepackage{enumitem}
\usepackage{cleveref}

\usepackage{scalerel,stackengine}
\stackMath
\newcommand\reallywidehat[1]{%
\savestack{\tmpbox}{\stretchto{%
  \scaleto{%
    \scalerel*[\widthof{\ensuremath{#1}}]{\kern-.6pt\bigwedge\kern-.6pt}%
    {\rule[-\textheight/2]{1ex}{\textheight}}
  }{\textheight}%
}{0.5ex}}%
\stackon[1pt]{#1}{\tmpbox}%
}

\theoremstyle{thmstyleone}%
\newtheorem{theorem}{Theorem}
\newtheorem{lemma}[theorem]{Lemma}
\newtheorem{proposition}[theorem]{Proposition}%
\newtheorem{assumption}{Assumption}
\theoremstyle{thmstyletwo}%
\newtheorem{remark}{Remark}%

\theoremstyle{thmstylethree}%

\newcommand{\ii}{\mathrm{i}}

\newcommand{\gao}{\bm{G}^{\alpha,\omega}}
\newcommand{\tgao}{\overline{\bm{G}}^{\alpha,\omega}}
\newcommand{\tgaop}[1]{\overline{G}^{\alpha,\omega}_{#1}}
\newcommand{\llm}{\mathcal{L}^{\lambda,\mu}}
\newcommand{\gpao}{\bm{G}_{+}^{\alpha,\omega}}

\begin{document}

\title[Subwavelength Resonance in Elastic Metascreen]{Mathematical Analysis of Subwavelength Resonance in Elastic Metascreen}


\author[1]{\fnm{Wei} \sur{Wu}}\email{wei\_wu@jlu.edu.cn}

\author*[2]{\fnm{Youzi} \sur{He}}\email{youzihe@smbu.edu.cn}


\affil[1]{\orgdiv{School of Mathematics}, \orgname{Jilin University}, \orgaddress{\street{2699 Qianjin Street}, \city{Changchun}, \postcode{130012}, \state{Jilin}, \country{China}}}

\affil*[2]{\orgdiv{MSU-BIT-SMBU Joint Research Center of Applied Mathematics}, \orgname{Shenzhen MSU-BIT University}, \orgaddress{\street{1 International University Park Road}, \city{Shenzhen}, \postcode{518172}, \state{Guangdong}, \country{China}}}

\affil*[2]{\orgdiv{School of Mathematics and Statistics}, \orgname{Beijing Institute of Technology}, \orgaddress{\city{Beijing}, \postcode{100081}, \state{Beijing}, \country{China}}}


\abstract{The aim of this paper is to provide a comprehensive and mathematically rigorous analysis on determining the existence of subwavelength resonance in elastic metascreen and resonance frequency calculation based on asymptotic analysis of quasi-periodic layer potential operators. An elastic metascreen is a thin sheet with subwavelength structures, which nevertheless has a significant effect on elastic wave propagation at specific frequencies. Periodic subwavelength elastic scatterers positioned on a reflective plane are considered in this paper. Firstly an explicit formula of quasi-periodic Green's function of Lam\'{e} system with Dirichlet boundary condition is derived for the first time. The subsequent discussion is twofold. In the first part where the shear modulus of scatterers is assumed to tend to infinity, the subwavelength resonance frequencies are given and approximated field inside inclusions and far-away from metascreen are calculated to demonstrate the dramatic change of scattered field due to subwavelength resonance. In the second part where the shear modulus of background is assumed to go to infinity, the absence of subwavelength resonance is proved. Without imposing conditions on the material parameters, the discussion in this paper provides the necessary condition for the occurrence of subwavelength resonance.}

\keywords{Subwavelength resonance, Quasi-periodic Green's function, Metascreen, Layer potential}



\maketitle

\section{Introduction}\label{sec1}

Metascreen is a special type of metamaterial, which refers to a thin sheet with designed periodic structures \cite{IEEE1}. These structures are typically manufactured to control the propagation, reflection, or scattering of waves by utilizing subwavelength features. One approach to design metascreen is to place identical microscopic inclusions in every unit of a periodic lattice. The size of these unit structures are smaller than the wavelength of the phenomena they affect. In acoustic metascreen, the unit size is usually of centimeter order, while for optic metascreen the size varies from millimeter order to nanometer order \cite{doi:10.1063/1.1343489,Shelby77,PhysRevLett.84.4184}. Elastic metascreen is a type of thin metamaterial engineered to manipulate elastic wave fields, in which the inclusions could be gas-filled, solid, or liquid micro-cavities. The physical nature of metascreen does not solely depend on the material of its unit structure. To a great extent it relies on the shape, size, alignment and orientation of its unit structure. By adjusting those properties, we could expect a remarkable change of corresponding metascreen in the absorption, enhancement or refraction to incident wave, of which we could take advantage to develop new materials possessing specific physical properties. In studies \cite{colombi2017enhanced, colombi2016seismic, kim2020elastic,zhu2014negative}, experimental evidence revealed the potential of these metascreens to achieve enhanced sensing, energy entrapment, and other marvelous results.

Recently, research has increasingly focused on the mathematical explanation of mechanism of exotic properties of different metamaterial and metascreen. For acoustic metascreen, \cite{ammari2017mathematical} gave a mathematical and numerical framework for the analysis and design of bubble meta-screens using layer potential and asymptotic analysis. \cite{feppon2023subwavelength} further analyzed the subwavelength resonance in finite acoustic one-dimensional periodic structure. Works like \cite{ammari2023mathematical, ammari2024fano} established mathematical explanation for electromagnetic metamaterial and metascreen. The studies concerning elastic meta-screen, however, remain relatively limited. \cite{ammari2006layer, chen2025analysis, ren2025subwavelength} provided results of high contrast asymptotics for the photonic crystals with priori knowledge on the contrast of $\lambda, \mu$ and $\rho$.

In this paper, we give a comprehensive study on subwavelength resonance in elastic metascreen. We consider an one-dimensional infinite metascreen composed of equally spaced identical elastic scatterers. The metascreen is placed horizontally in two-dimensional space, slightly above a total-reflective plane. The scatterers and background medium possesses different Lam\'{e} constants $\lambda,\mu$ and density $\rho$. As the basis of all subsequent calculations in this paper, we firstly derive the explicit formula of quasi-periodic Green's function of Lam\'{e} system. As shown in existing researches, the occurrence of subwavelength resonance must be resulted from high contrast of physical parameters. Hence we consider two different limiting cases here, the shear modulus of inclusions tending to infinity and the shear modulus of the background going to infinity, respectively. We do not impose any further conditions on the other Lam\'{e} constants and densities. In the first case where the shear modulus of inclusions goes to infinity, we will show the requirements that the remaining physical constants must satisfy to ensure the appearance of subwavelength resonance. We will demonstrate the approximation of fields inside the scatterer and far-away from metascreen, from which we can observe that both near field and far field dramatically alters when subwavelength resonance happens. In the second case we will prove that the subwavelength resonance will not emerge at all events. Inspired by \cite{ammari2017mathematical, HabibARMA, li2023mathematicaltheorydipolarresonances}, we will use asymptotic spectral analysis on layer potential operators along with general Rouch\'{e}'s theorem. The equation describing elastic wave in $\mathbb{R}^2$ will be reformulated as an integral equation system comprised of layer potential operators. Finding subwavelength resonance then equivalents to finding characteristic value of the integral system. We apply asymptotic expansion to calculate an approximated integral system, whose characteristic value is much easier to find. This paper provides a complete answer to what kind of elastic metamaterials with high shear modulus contrast can exhibit subwavelength resonance. The same method could be applied to the cases of high density contrast and high bulk modulus contrast. The discussions here could be generalized to the one-dimensional and two-dimensional elastic metascreens in three-dimensional space without much difficulty.

The remainder of the paper is organized as follows. In Section 2, we present the mathematical formulation of elastic metascreen and the incident elastic wave. In Section 3, we give the explicit expression of quasi-periodic Green's function of Lam\'{e} system with Dirichlet boundary condition and introduce layer potential operators with some important properties for subsequent discussions. With the help of layer potential operators, we transform the problem into a system of integral equations. In Section 4, we deduce the existence of sub-wavelength resonances when the shear modulus of inclusions tending to infinity. We also calculate the near-field and far-field approximation of elastic field, from which we could tell that the field inside scatterer and far-away from metascreen drastically varies on the occasion of subwavelength resonance. The absence of subwavelength resonance for the case that the shear modulus of the background going to infinity is presented in Section 5.

\section{Statement of the problem}\label{sec2}

We consider the one-dimensional quasi-periodic elastic metascreen in two-dimensional plane aligned along $x_1$-axis. Let $\partial \mathbb{R}^2_+:=\{(x_1,x_2)\in\mathbb{R}^2, x_2=0 \}$ represent the boundary of half plane and $\mathbb{R}^2_+:=\{\bm{x}\in\mathbb{R}^2, x_2>0 \}$ be the upper half-space. Let $D$ be a connected domain with the Lipschitz boundary lying inside the periodic cell $S:=[-1/2,1/2)\times\mathbb{R}_+$. In the following, $D$ is regarded as a periodic repeated unit. In our model, the identical inclusions are periodically aligned along $x_1$-axis with periodicity 1. We call the array of all periodic aligned inclusions $\Omega:=\bigcup_{n\in\mathbb{Z}}(D+n)$.
\begin{figure}[H]
\centering
\includegraphics[width=0.7\textwidth]{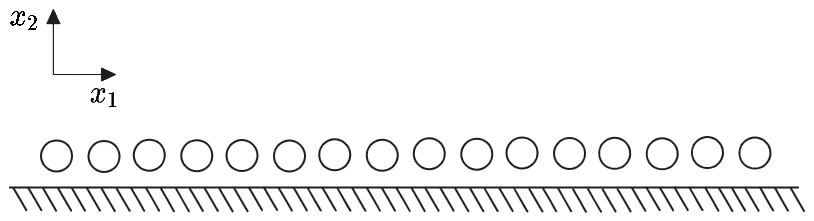}
\caption{Schematic figure of elastic metascreen.}
\end{figure}
As we know, the elastic nature of medium is described by its Lam\'{e} constants and density. We assume here that the metascreen $\Omega$ and background medium are homogeneous with different physical constants. The Lam\'{e} constants inside $\Omega$ are $\tilde{\lambda}, \tilde{\mu}$, and the Lam\'{e} constants of the background elastic medium are $\lambda$ and $\mu$. The medium density of background and inclusions are $\rho$ and $\widetilde{\rho}$, correspondingly. The elastic scattering problem can be modeled by the following system
\begin{equation}\label{eq:elas_sys}
	\left\{
	\begin{array}{ll}
		(\mathcal{L}^{\lambda, \mu} +\rho\omega^2) \bm{u} = \bm{0},  &   {\rm{in }}\ \mathbb{R}^2\backslash\overline{\Omega}, \medskip \\
		\displaystyle  (\mathcal{L}^{\tilde{\lambda}, \tilde{\mu}} +\widetilde{\rho}\omega^2) \bm{u} = \bm{0},  &  {\rm{in }}\  \Omega, \medskip \\
		\displaystyle  \bm{u}|_+ - \bm{u}|_- = \bm{0},   & \rm{on }\  \partial \Omega, \medskip \\
                 \displaystyle \frac{\partial\bm{u}}{\partial\bm{\nu} }\bigg|_+ - \frac{\partial\bm{u}}{\partial\tilde{\bm{\nu}} }\bigg|_-  = \bm{0},   &\rm{on }\  \partial \Omega, \medskip \\
                 \displaystyle \bm{u} = \bm{0},  &\rm{on }\   \partial \mathbb{R}^2_+, \medskip \\
		\displaystyle \bm{u}(x_1+1,x_2) = e^{\ii\alpha}\bm{u}(x_1,x_2).
	\end{array}
	\right.
\end{equation}
Here constant $\alpha$ denotes the phase change between neighboring cells. $\mathcal{L}^{\lambda, \mu}$ is the Lam\'{e} operator corresponding to the Lam\'{e} constants $\lambda, \mu$ defined by
\begin{equation*}
\mathcal{L}^{\lambda, \mu}\bm{u} := \mu\Delta \bm{u} +(\lambda+ \mu)\nabla\nabla\cdot \bm{u},
\end{equation*}
and the corresponding conormal derivative $\partial\bm{u}/\partial\bm{\nu} $ on $\partial\Omega$ is defined to be
\begin{equation*}
\frac{\partial\bm{u}}{\partial\bm{\nu} }:=\lambda(\nabla\cdot \bm{u})\bm{\nu} +\mu(\nabla\bm{u}+\nabla\bm{u}^{\top} )\bm{\nu},
\end{equation*}
where $\bm{\nu}$ is the unit outward normal to $\partial\Omega$ and the superscript $\top$ denotes the transpose of a matrix. The notation $\partial\bm{u}/\partial\bm{\widetilde{\nu}}$ denotes the conormal derivative with Lam\'{e} constants $\widetilde{\lambda}, \widetilde{\mu}$.

Let $\bm{u}^{\rm{in}}=\bm{u}^{\rm{in}}(\bm{x})$ be the incident wave. We consider the metascreen is shined by planar P-wave, which satisfies
\begin{equation*}
\bm{u}^P(\bm{x})= e^{\ii\bm{k}_P\cdot\bm{x}}\bm{\theta} = {\rm{e}}^{{\rm{i}}\omega\sqrt{\rho/(\lambda+2\mu)} \bm{\theta}\cdot\bm{x}}\bm{\theta}
\end{equation*}
where $\bm{k}_P:= \omega\sqrt{\rho/(\lambda+2\mu)} \bm{\theta} = (k_{P,1}, k_{P,2})$ is the wave vector of planar P-wave. we denote the corresponding wave speed as $c_p = \sqrt{(\lambda+2\mu)/\rho}$ and wavenumber as $k_p = \omega\sqrt{\rho/(\lambda+2\mu)}$. The wave vector of planar S-wave is $\bm{k}_S:=\omega\sqrt{\rho/\mu}\bm{\theta}^\perp$, with wave speed $c_s=\sqrt{\mu/\rho}$ and wavenumber $k_s=\omega\sqrt{\rho/\mu}$. It is worth mentioning that although we only consider P-wave in this paper, the same steps could be applied to planar S-wave. Since in \eqref{eq:elas_sys} we imposed Dirichlet boundary condition on $\mathbb{R}^2_+$, the incident wave is reflected by $\mathbb{R}^2_+$, so the actual incident wave applied onto metascreen is
\begin{equation}\label{eq:upp}
\bm{u}^P_+(\bm{x}):= \bm{u}^P(x_1,x_2)-\bm{u}^P(x_1,-x_2) = 2\ii{\rm{e}}^{\ii\omega\sqrt{\frac{\rho}{\lambda+2\mu}}\theta_1x_1}\sin(\omega\sqrt{\frac{\rho}{\lambda+2\mu}}\theta_2x_2)\bm{\theta}.
\end{equation}

Here we calculate the expression of $\partial\bm{u}_+^P/\partial\bm{\nu}$ for later use. From definition of boundary contraction we have
\begin{equation*}
\left(\frac{\partial\bm{u}_+^P}{\partial\bm{\nu} }\right)_k = \lambda(\nabla\cdot\bm{u}_+^P)\nu_k + \mu\sum_{j=1}^{2}  \Big(\frac{\partial u_{+,k}^P}{\partial x_j }\nu_j + \frac{\partial u_{+,j}^P}{\partial x_k }\nu_j\Big), k=1,2.
\end{equation*}
The explicit expression of the divergence term is
\begin{equation}\label{eq:upnu1}
\begin{aligned}
\nabla\cdot\bm{u}_+^P(\bm{x})=& 2{\rm{i}}{\rm{e}}^{{\rm{i}}k_{P,1}x_1}\text{sin}(k_{P,2}x_2){\rm{i}}k_{P,1}\theta_1 + 2{\rm{i}}{\rm{e}}^{{\rm{i}}k_{P,1}x_1}k_{P,2}\text{cos}(k_{P,2}x_2)\theta_2\\
=&2\omega\sqrt{\frac{\rho}{\lambda+2\mu}}e^{\ii\omega\sqrt{\frac{\rho}{\lambda+2\mu}}\theta_1x_1}\left(-\theta_1^2\sin(\omega\sqrt{\frac{\rho}{\lambda+2\mu}}\theta_2x_2) + \ii\theta_2^2\cos(\omega\sqrt{\frac{\rho}{\lambda+2\mu}}\theta_2x_2)\right)
\end{aligned}
\end{equation}
For the second part, we calculate that
\begin{equation*}
\frac{\partial u^{P}_{+,k}}{\partial x_j}=\frac{\partial}{\partial x_j}\big(2{\rm{i}}{\rm{e}}^{{\rm{i}}k_{P,1}x_1}\text{sin}(k_{P,2}x_2)\theta_k\big)=
\left\{
	\begin{array}{ll}
-2k_{P,1}{\rm{e}}^{{\rm{i}}k_{P,1}x_1}\sin(k_{P,2}x_2)\theta_k,& j=1,\\
2\ii k_{P,2}{\rm{e}}^{{\rm{i}}k_{P,1}x_1}\text{cos}(k_{P,2}x_2)\theta_k,& j=2.
\end{array}
	\right.
\end{equation*}
So when $k=1$, by direct calculation we can obtain
\begin{equation}\label{eq:upnu21}
\begin{aligned}
\sum_{j=1}^{2}  \Big(\frac{\partial u_{+,1}^P}{\partial x_j } + \frac{\partial u_{+,j}^P}{\partial x_1 }\Big)\nu_j = &  \Big(\frac{\partial u_{+,1}^P}{\partial x_1 } + \frac{\partial u_{+,1}^P}{\partial x_1 }\Big)\nu_1 + \Big(\frac{\partial u_{+,1}^P}{\partial x_2 } + \frac{\partial u_{+,2}^P}{\partial x_1 }\Big)\nu_2\\
=&\big(-2k_{P,1}{\rm{e}}^{{\rm{i}}k_{P,1}x_1}\text{sin}(k_{P,2}x_2)\theta_1 -2k_{P,1}{\rm{e}}^{{\rm{i}}k_{P,1}x_1}\text{sin}(k_{P,2}x_2)\theta_1  \big)\nu_1\\
&+\big(2{\rm{i}}k_{P,2}{\rm{e}}^{{\rm{i}}k_{P,1}x_1}\text{cos}(k_{P,2}x_2)\theta_1 -2k_{P,1}{\rm{e}}^{{\rm{i}}k_{P,1}x_1}\text{sin}(k_{P,2}x_2)\theta_2\big)\nu_2\\
=&-4k_{P,1}{\rm{e}}^{{\rm{i}}k_{P,1}x_1}\text{sin}(k_{P,2}x_2)\theta_1\nu_1\\
&+2e^{\ii k_{P,1}x_1}(\ii k_{P,2}\cos(k_{P,2}x_2)\theta_1 - k_{P,1}\sin(k_{P,2}x_2)\theta_2)\nu_2.
\end{aligned}
\end{equation}
When $k= 2$, one has
\begin{equation}\label{eq:upnu22}
\begin{aligned}
&\sum_{j=1}^{2} \Big(\frac{\partial u_{+,2}^P}{\partial x_j } + \frac{\partial u_{+,j}^P}{\partial x_2 }\Big)\nu_j\\
=& 2e^{\ii k_{P,1}x_1}(-k_{P,1}\sin(k_{P,2}x_2)\theta_2 + \ii k_{P,2}\cos(k_{P,2}x_2)\theta_1)\nu_1 + 4\ii k_{P,2}e^{ik_{P,1}x_1}\cos(k_{P,2}x_2)\theta_2\nu_2.
\end{aligned}
\end{equation}
Combining \eqref{eq:upnu1}, \eqref{eq:upnu21} and \eqref{eq:upnu22} we can obtain the explicit formula for $\partial\bm{u}_+^P/\partial\bm{\nu}$.


\section{Preliminaries}\label{sec3}

In this section, we aim at defining the quasi-periodic Green's function of Lam\'{e} system with Dirichlet boundary condition on $\mathbb{R}^2_+:=\{(x_1,0):x_1\in\mathbb{R}\}$. We derive the explicit expression of the Green's function and separate the evanescent and propagating parts. Based on the aforementioned Green's function we define quasi-periodic layer potential operators and introduce some important properties, which will be very useful in subsequent discussions.

Throughout the whole paper, we consider unitary Fourier transform, i.e.
$$
    \hat{f}(\xi):=\frac{1}{(2\pi)^{n/2}}\int_{\mathbb{R}^n} f(x) e^{-\ii x\cdot\xi}{\rm{d}} x
$$
for $x,\xi\in\mathbb{R}^n$.

\subsection{Quasi-periodic-Dirichlet Green's function}\label{sec3-1}

Firstly we give the definition and properties of quasi-periodic-Dirichlet Green's function for the Lam\'{e} system. We first introduce the quasi-periodic Green's function $\bm{G}^{\alpha,\omega}$ satisfying
\begin{equation}\label{eq:quasi-p-Green}
\begin{split}
(\mathcal{L}^{\lambda, \mu} +\rho\omega^2\bm{I})\bm{G}^{\alpha,\omega}(\bm{x},\bm{y})&=\sum_{n \in \mathbb{Z}} \delta(x_1 - y_1 - n) \delta(x_2 - y_2) {\rm{e}}^{{\rm{i}} n\alpha}\bm{I}\\
&=\delta(x_2 - y_2)\left(\sum_{n \in \mathbb{Z}} \delta(x_1 - y_1 - n){\rm{e}}^{{\rm{i}} n\alpha} \right)\bm{I}.
\end{split}
\end{equation}

Here $\bm{I}$ is $2\times 2$ identity matrix, and $\bm{G}^{\alpha,\omega}(\bm{x},\bm{y})$ is $2\times 2$ matrix-valued function. It is straightforward that $\bm{G}^{\alpha,\omega}(\bm{x},\bm{y}) = \bm{G}^{\alpha,\omega}(\bm{x}-\bm{y})$. Define $\overline{\bm{G}}^{\alpha,\omega}(\bm{x}):= e^{-\ii \alpha x_1}\bm{G}^{\alpha,\omega}(\bm{x})$. $\tgao(\bm{x})$ is then a periodic function of $x_1$ with periodicity 1. Substituting $e^{\ii \alpha x_1}\overline{\bm{G}}$ into \eqref{eq:quasi-p-Green} we get the following results.
    \begin{equation}
        \begin{aligned}
            &e^{-\ii\alpha x_1}(\llm(e^{\ii\alpha x_1} \tgao))_{1j} = \\&\mu(-\alpha^2+2\ii\alpha\partial_1 + \Delta)\tgaop{1j}+(\lambda+\mu)(-\alpha^2+2\ii\alpha\partial_1+\partial_{11})\tgaop{1j} + (\lambda+\mu)(\ii\alpha\partial_2+\partial_{12})\tgaop{2j}, \\
            &e^{-\ii\alpha x_1}(\llm(e^{\ii\alpha x_1}\tgao))_{2j} = \\
            &\mu(-\alpha^2+2\ii\alpha\partial_1+\Delta)\tgaop{2j}+(\lambda+\mu)(\ii\alpha\partial_2+\partial_{12})\tgaop{1j} +(\lambda+\mu)\partial_{22}\tgaop{2j}.
        \end{aligned}
    \end{equation}

    Here $j=1,2$, and $\partial_if$ is short for $\frac{\partial f}{\partial x_i}$. Define
$$
    \overline{\mathcal{L}}^{\lambda,\mu}:=\begin{pmatrix}
        \mu\partial_{22}+(\lambda+2\mu)(-\alpha^2+2\ii\alpha\partial_1+\partial_{11}) & (\lambda+\mu)(\ii\alpha\partial_2+\partial_{12}) \\
        (\lambda+\mu)(\ii\alpha\partial_2+\partial_{12}) & (\lambda+2\mu)\partial_{22}+\mu(-\alpha^2+2\ii\alpha\partial_1+\partial_{11})
    \end{pmatrix},
$$
then
\begin{equation}\label{eq:tgao1}
    (\overline{\mathcal{L}}^{\lambda,\mu}+\rho\omega^2\bm{I})\tgao(\bm{x}) = \sum\limits_{n\in\mathbb{Z}}\delta(x_1-n)\delta(x_2)\bm{I}.
\end{equation}
Consider the Fourier expansion of $\tgao(\bm{x})$ with respect to $x_1$ in the form of
\begin{equation}\label{eq:tgaopij}
    \tgaop{ij}(x_1,x_2) = \sum\limits_{l\in 2\pi\mathbb{Z}} c_{l,ij}(x_2)e^{\ii lx_1}.
\end{equation}
To get rid of $x_1$ dependency in the equation, we take advantage of Poisson summation formula
$$
    \sum\limits_{n\in\mathbb{Z}}\delta(x_1-n) = \frac{1}{2\pi}\sum\limits_{l\in 2\pi\mathbb{Z}}e^{\ii lx_1}
$$
Plugging Fourier expansion and Poisson summation formula into \eqref{eq:tgao1} leads to an ordinary differential equation system of $\bm{c}_l(x_2):=(c_{l,ij}(x_2))_{2\times 2}$.The expressions are shown below.

\begin{lemma}\label{eq:cl}
    $\bm{c}_l$ satisfies
    \begin{equation}\label{eq:clode}
        \begin{aligned}
            (-(\lambda+2\mu)(\alpha+l)^2+\rho\omega^2)c_{l,1j} + \mu c_{l,1j}'' + \ii(\lambda+\mu)(\alpha+l)c_{l,2j}' &= \frac{1}{2\pi}\delta_0\delta_{1j} \\
            (-\mu(\alpha+l)^2+\rho\omega^2)c_{l,2j} + (\lambda+2\mu)c_{l,2j}'' + \ii(\lambda+\mu)(\alpha+l)c_{l,1j}' &= \frac{1}{2\pi}\delta_0\delta_{2j}
        \end{aligned}
    \end{equation}
    for $j=1,2$.
\end{lemma}

Notice that $c_{l,11}$, $c_{l,21}$ (as well as $c_{l,12}$, $c_{l,22}$) are paired in the equation system. For brevity we only demonstrate the calculation details of $c_{l,11}$ and $c_{l,21}$, and directly give the final expressions of $c_{l,12}$ and $c_{l,22}$.

Applying Fourier transform to \eqref{eq:clode}, we get
\begin{equation*}\label{eq:clft}
    \begin{aligned}
        (\rho\omega^2-\mu(\alpha+l)^2-(\lambda+2\mu)\xi^2)\hat{c}_{l,21} - (\lambda+\mu)(\alpha+l)\xi\hat{c}_{l,11} &= 0 \\
        (\rho\omega^2-(\lambda+2\mu)(\alpha+l)^2-\mu\xi^2)\hat{c}_{l,11} - (\lambda+\mu)(\alpha+l)\xi\hat{c}_{l,21} &= (2\pi)^{-\frac{3}{2}}.
    \end{aligned}
\end{equation*}
Define
\begin{align*}
D &= (\rho\omega^2-(\lambda+2\mu)(\alpha+l)^2-\mu\xi^2)(\rho\omega^2-\mu(\alpha+l)^2-(\lambda+2\mu)\xi^2)-(\lambda+\mu)^2(\alpha+l)^2\xi^2. \\
&= \mu(\lambda+2\mu)(\xi^2+(\alpha+l)^2-\frac{\rho\omega^2}{\lambda+2\mu})(\xi^2+(\alpha+l)^2-\frac{\rho\omega^2}{\mu}).
\end{align*}

Calculation shows that
\begin{equation*}
\begin{aligned}
\hat{c}_{l,11}(\xi) &= (2\pi)^{-3/2}  \frac{\rho\omega^2-\mu(\alpha+l)^2-(\lambda+2\mu)\xi^2}{D}, \\
\hat{c}_{l,21}(\xi)
&= (2\pi)^{-3/2}\frac{(\lambda+\mu)(\alpha+l)\xi}{D}.
\end{aligned}
\end{equation*}

Here, $\alpha$ is exactly the first element of wave vector $\bm{k}$, so $\alpha<k_p$ and $\alpha<k_s$.

We define
\begin{align*}
l_{p,+}:= k_p-\alpha, && l_{p,-}:=-k_p-\alpha, \\
l_{s,+}:= k_s-\alpha, && l_{s,-}:=-k_s-\alpha.
\end{align*}
It is obvious that $l_{p,+}<l_{s,+}$ and $l_{p,-}>l_{s,-}$ since $k_p<k_s$. Define $\gamma_{l,s} := \sqrt{|(\alpha+l)^2-\frac{\rho\omega^2}{\mu}|}, \gamma_{l,p} := \sqrt{|(\alpha+l)^2-\frac{\rho\omega^2}{\lambda+2\mu}|}$. The following calculation is divided into three different cases:
\begin{itemize}
\item \textbf{Case I: $l\geq l_{s,+}$ or $l\leq l_{s,-}$.}

In this case $(\alpha+l)^2\geq k_s^2>k_p^2$. To calculate the inverse Fourier transform of $\hat{c}_{l,ij}$, we need the following result.
\begin{lemma}
    \begin{equation*}
    \begin{aligned}
        \int_{-\infty}^{+\infty} \frac{e^{\ii\xi x}}{\xi^2+a^2}\mathrm{d}\xi &= \pi\frac{e^{-|ax|}}{a}, \\
        \int_{-\infty}^{+\infty} \frac{\xi e^{\ii\xi x}}{\xi^2+a^2}\mathrm{d}\xi &= \ii\pi\mathrm{sgn}(x)e^{-|ax|}.
    \end{aligned}
    \end{equation*}
\end{lemma}

With the help of lemma above we conclude that
\begin{align}\label{eq:cl1121}
    c_{l,11}(x) &= \frac{\lambda+\mu}{4\pi\mu(\lambda+2\mu)(\gamma_{l,p}^2-\gamma_{l,s}^2)}\left(\gamma_{l,s}e^{-\gamma_{l,s}|x|} - \frac{(\alpha+l)^2}{\gamma_{l,p}}e^{-\gamma_{l,p}|x|}\right),\\
    c_{l,21}(x) &= \frac{\ii(\lambda+\mu)(\alpha+l)}{4\pi\mu(\lambda+2\mu)(\gamma_{l,p}^2-\gamma_{l,s}^2)}\mathrm{sgn}(x)(e^{-\gamma_{l,s}|x|} - e^{-\gamma_{l,p}|x|}).
\end{align}

Similarly, we can solve the other two Fourier coefficients $c_{l,12}, c_{l,22}$ satisfying
\begin{equation*}
    \begin{aligned}
        (\rho\omega^2-\mu(\alpha+l)^2-(\lambda+2\mu)\xi^2)\hat{c}_{l,22} - (\lambda+\mu)(\alpha+l)\xi\hat{c}_{l,12} &= (2\pi)^{-\frac{3}{2}}, \\
        (\rho\omega^2-(\lambda+2\mu)(\alpha+l)^2-\mu\xi^2)\hat{c}_{l,12} - (\lambda+\mu)(\alpha+l)\xi\hat{c}_{l,22} &= 0.
    \end{aligned}
\end{equation*}

as follows
\begin{align}\label{eq:cl1222}
c_{l,12}(x) &= \frac{\ii(\lambda+\mu)(\alpha+l)}{4\pi\mu(\lambda+2\mu)(\gamma_{l,p}^2-\gamma_{l,s}^2)}\mathrm{sgn}(x)(e^{-\gamma_{l,s}|x|}-e^{-\gamma_{l,p}|x|}), \\
c_{l,22}(x) &= \frac{\lambda+\mu}{4\pi\mu(\lambda+2\mu)(\gamma_{l,p}^2-\gamma_{l,s}^2)}\left(\gamma_{l,p}e^{-\gamma_{l,p}|x|} - \frac{(\alpha+l)^2}{\gamma_{l,s}}e^{-\gamma_{l,s}|x|}\right).
\end{align}

We can see that all $c_{l,ij}$ are decaying exponentially with the increasing $|x|$. This implies that the corresponding part in Green's function generate an evanescent wave. As a result, all modes with $l\geq l_{s,+}$ or $l\leq l_{s,-}$ are evanescent modes.

\item \textbf{Case II: $l_{p,+}\leq l<l_{s,+}$ or $l_{s,-}< l \leq l_{p,-}$.}

In this case $k_s^2>(\alpha+l)^2\geq k_p^2$. To calculate the inverse Fourier transform of $c_{l,ij}$ we need two additional results from Sommerfeld integrals.
\begin{lemma}\label{lem:feymann}
    With Feymann prescription, the following identities hold.
    \begin{align*}
       \int_{-\infty}^{+\infty} \frac{e^{ix\xi}}{\xi^2-k^2}\mathrm{d}\xi &= \frac{\ii\pi}{k}e^{\ii k|x|}, \\
       \int_{-\infty}^{+\infty} \frac{\xi e^{ix\xi}}{\xi^2-k^2}\mathrm{d}\xi &= \ii\pi\mathrm{sgn}(x)e^{\ii k|x|}.
    \end{align*}
\end{lemma}
Then with a bit calculation we conclude that
\begin{equation}\label{eq:lplls}
\begin{aligned}
    c_{l,11}(x) &= \frac{-(\lambda+\mu)}{4\pi(\lambda+2\mu)\mu(\gamma_{l,s}^2+\gamma_{l,p}^2)}\left(\frac{(\alpha+l)^2}{\gamma_{l,p}}e^{-\gamma_{l,p}|x|} + \ii\gamma_{l,s}e^{\ii\gamma_{l,s}|x|}\right), \\
    c_{l,21}(x) &= c_{l,12}(x) =  \frac{-\ii(\alpha+l)(\lambda+\mu)\mathrm{sgn}(x)}{4\pi(\lambda+2\mu)\mu(\gamma_{l,s}^2+\gamma_{l,p}^2)}(e^{-\gamma_{l,p}|x|} - e^{\ii\gamma_{l,s}|x|}), \\
    c_{l,22}(x) &= \frac{-(\lambda+\mu)}{4\pi(\lambda+2\mu)\mu(\gamma_{l,s}^2+\gamma_{l,p}^2)}\left(-\gamma_{l,p}e^{-\gamma_{l,p}|x|} + \frac{\ii(\alpha+l)^2}{\gamma_{l,s}}e^{\ii\gamma_{l,s}|x|}\right).
\end{aligned}
\end{equation}

We can see that the p-wave part in the Green's function is exponentially decreasing, while the s-wave part propagates. Therefore when $l_{p,+}\leq l<l_{s,+}$ or $l_{s,-}< l \leq l_{p,-}$, only s-wave contributes to the propagating solution, and p-wave stays evanescent.

\item \textbf{Case III: $l_{p,-}<l<l_{p,+}$.}
When $l_{p,-}<l<l_{p,+}$, the term $\rho\omega^2-\mu\alpha^2>0$ and $\rho\omega^2-(\lambda+2\mu)\alpha^2>0$. With the help of Lemma \ref{lem:feymann}
we can get that
\begin{equation}\label{eq:c0}
\begin{aligned}
    c_{l,11}(x) &= \frac{\ii(\lambda+\mu)}{4\pi\mu(\lambda+2\mu)(\gamma_{l,p}^2-\gamma_{l,s}^2)}\left(\gamma_{l,s}e^{i\gamma_{l,s}|x|} + \frac{(\alpha+l)^2}{\gamma_{l,p}}e^{i\gamma_{l,p}|x|}\right), \\
    c_{l,12}(x) &= c_{l,21}(x) = \frac{\ii(\lambda+\mu)(\alpha+l)}{4\pi\mu(\lambda+2\mu)(\gamma_{l,p}^2-\gamma_{l,s}^2)}\mathrm{sgn}(x)(e^{\ii \gamma_{l,p}|x|} - e^{\ii \gamma_{l,s}|x|}), \\
    c_{l,22}(x) &=
    \frac{\ii(\lambda+\mu)}{4\pi\mu(\lambda+2\mu)(\gamma_{l,p}^2-\gamma_{l,s}^2)}\left(\gamma_{l,p}e^{\ii\gamma_{l,p}|x|} + \frac{(\alpha+l)^2}{\gamma_{l,s}}e^{\ii\gamma_{l,s}|x|}\right).
\end{aligned}
\end{equation}
This leads to a propagating solution for both p-wave and s-wave instead of an evanescent one.
\end{itemize}

To conclude, with \eqref{eq:cl1121}, \eqref{eq:cl1222}, \eqref{eq:lplls}, \eqref{eq:c0}, we achieved the full expression of $\tgao(\bm{x})$. The desired quasi-periodic Green's function satisfying \eqref{eq:quasi-p-Green} could be calculated by
$$
\gao(\bm{x}) = e^{\ii\alpha x_1}\tgao(\bm{x}) = \sum\limits_{l\in 2\pi\mathbb{Z}}e^{\ii (l+\alpha)x_1}\begin{pmatrix}
c_{l,11}(x_2) & c_{l,12}(x_2) \\ c_{l,21}(x_2) & c_{l,22}(x_2)
\end{pmatrix} = \sum\limits_{l\in 2\pi\mathbb{Z}}e^{\ii (l+\alpha)x_1}\bm{c}_l(x_2).
$$
Since in our discussion, there is a Dirichlet boundary condition imposed on $\partial\mathbb{R}^2_+$, to facilitate our discussion we define $\gpao(\bm{x},\bm{y}):= \gao((x_1-y_1,x_2-y_2))-\gao((x_1-y_1,x_2+y_2))$. With this definition $\gpao(\bm{x},\bm{y})$ takes the following form
\begin{equation}\label{eq:quasiperiodicgreen}
    \gpao(\bm{x},\bm{y}) = \sum\limits_{l\in 2\pi\mathbb{Z}}(\bm{c}_l(x_2-y_2)-\bm{c}_l(x_2+y_2))e^{\ii(\alpha+l)(x_1-y_1)}
\end{equation}
and $\gpao(\bm{x},\bm{y}) = 0$ when $\bm{x}\in\mathbb{R}_+^2$ and $\bm{y}\in\partial\mathbb{R}_+^2$. Notice that unlike $\gao$, $\gpao$ is no longer translation invariant, i.e. $\gpao(\bm{x},\bm{y})\neq \gpao(\bm{x}-\bm{y},0)$.

\subsection{Quasi-periodic-Dirichlet layer potentials}\label{sec3-2}

Now, we are ready to introduce the quasi-periodic-Dirichlet layer potential operators. For $\bm{\psi}\in (L^2(\partial D))^2$, we define the following layer potential operators with quasi-periodic Green's function defined in \eqref{eq:quasiperiodicgreen}.
\begin{equation}\label{eq:layer-potential-D}
\begin{split}
\mathcal{S}^{\alpha,\omega}_+[\bm{\psi}](\bm{x}):=&\int_{\partial D} \bm{G}_{+}^{\alpha,\omega} (\bm{x},\bm{y})\bm{\psi}(\bm{y}){\rm{d}}\sigma(\bm{y}), \ \bm{x}\in \mathbb{R}_+^2,\\
\mathcal{K}^{\alpha,\omega}_+[\bm{\psi}](\bm{x}):=&\mathrm{p.v.}\int_{\partial D} \frac{\partial\bm{G}_{+}^{\alpha,\omega} }{\partial \bm{\nu_x}}(\bm{x},\bm{y})\bm{\psi}(\bm{y}){\rm{d}}\sigma(\bm{y}),\  \bm{x}\in \partial D,\\
\widetilde{\mathcal{S}}^{\alpha,\omega}_+[\bm{\psi}](\bm{x}):=&\int_{\partial D} \widetilde{\bm{G}}_{+}^{\alpha,\omega} (\bm{x},\bm{y})\bm{\psi}(\bm{y}){\rm{d}}\sigma(\bm{y}), \ \bm{x}\in \mathbb{R}^2_+,\\
\widetilde{\mathcal{K}}^{\alpha,\omega}_+[\bm{\psi}](\bm{x}):=&\mathrm{p.v.}\int_{\partial D} \frac{\partial\widetilde{\bm{G}}_{+}^{\alpha,\omega} }{\partial \bm{\nu_x}}(\bm{x},\bm{y})\bm{\psi}(\bm{y}){\rm{d}}\sigma(\bm{y}),\  \bm{x}\in \partial D.\\
\end{split}
\end{equation}
Here $\widetilde{\bm{G}}^{\alpha,\omega}_+$ denotes the $\alpha$-quasi-periodic Dirichlet Green's function corresponding to Lam\'{e} system with Lam\'{e} constants $\tilde{\lambda}, \tilde{\mu}$. $\partial/\partial{\bm{\nu}_x}$ means taking normal derivative with respect to $\bm{x}$.

The following properties of $\mathcal{S}^{\alpha,\omega}_+$ and $\mathcal{K}^{\alpha,\omega}_+$ could be proved with similar methods as in \cite{HabibElastic}.
\begin{proposition}\label{prop:layerpotentialproperty}
Single layer potential $\mathcal{S}^{\alpha,\omega}_+[\bm{\psi}]$ is continuous across $\partial D$, i.e.
$$
    \left.\mathcal{S}^{\alpha,\omega}_+[\bm{\psi}]\right\vert_- = \left.\mathcal{S}^{\alpha,\omega}_+[\bm{\psi}]\right\vert_+ \quad \mathrm{on}~\partial D.
$$
When regarding as operator from $L^2(\partial D)$ to $L^2(\partial D)$, the $L^2$ adjoint of $\mathcal{K}^{\alpha,\omega}_+$ is
$$
    (\mathcal{K}^{\alpha,\omega}_+)^*[\bm{\psi}](\bm{x}) := \mathrm{p.v.}\int_{\partial D} \frac{\partial\bm{G}_{+}^{\alpha,\omega} }{\partial \bm{\nu_y}}(\bm{x},\bm{y})\bm{\psi}(\bm{y}){\rm{d}}\sigma(\bm{y}),\  \bm{x}\in \partial D.
$$
The normal derivative of the single layer potential satisfies the following jump relation
\begin{equation*}
\left.\frac{\partial}{\partial\bm{\nu}} \mathcal{S}_+^{\alpha,\omega}[\bm{\psi}]\right\vert_{\pm}(\bm{x})=\Big(\pm\frac{1}{2}\bm{I}+(\mathcal{K}_+^{\alpha,\omega})^{\ast} \Big)\bm{\psi}(\bm{x})\quad {\rm{a.e.}}\ \bm{x}\in \partial D.
\end{equation*}
\end{proposition}
With the help of \eqref{eq:layer-potential-D} and Proposition \ref{prop:layerpotentialproperty}, we could write the solution to problem \eqref{eq:elas_sys} in the following integral representation formula
\begin{equation}\label{eq:int_repre}
\bm{u}(\bm{x}):=\left\{
	\begin{array}{ll}
		\bm{u}^{{\rm{in}}}_+(\bm{x})+\mathcal{S}^{\alpha,\omega}_+[\bm{\psi}](\bm{x}),  &  \bm{x} \in  \mathbb{R}^2\setminus\bar{\Omega}, \medskip \\
                 \displaystyle  \widetilde{\mathcal{S}}^{\alpha,\omega}_+[\bm{\phi}](\bm{x}) ,   &\bm{x} \in  \Omega ,
	\end{array}
	\right.
\end{equation}
The boundary potentials $\bm{\phi}, \bm{\psi}\in (L^2(\partial D))^2$ are determined by
\begin{equation*}\label{eq:A_F}
\mathcal{A}\left(
\begin{array}{c}
     \bm{\phi} \\
      \bm{\psi}
\end{array}
\right)=\bm{F},
\end{equation*}
where
\begin{equation*}
\mathcal{A}:=\left(
\begin{array}{cc}
\widetilde{\mathcal{S}}^{\alpha,\omega}_+ & -\mathcal{S}^{\alpha,\omega}_+\\
-\frac{1}{2}\bm{I}+(\widetilde{\mathcal{K}}_+^{\alpha,\omega})^{\ast} &-\big(\frac{1}{2}\bm{I}+(\mathcal{K}_+^{\alpha,\omega})^{\ast} \big)
\end{array}
\right)
\end{equation*}
and
\begin{equation}\label{eq:Fdefinition}
\bm{F}:=\left(
\begin{array}{c}
    \bm{u}^{{\rm{in}}}_+  \\
      \frac{\partial}{\partial\bm{\nu}}\bm{u}^{{\rm{in}}}_+
\end{array}
\right).
\end{equation}

\section{Elastic metascreen with $\widetilde{\mu}\rightarrow +\infty$}

In this section, we will study the behavior of elastic metascreen when the shear modulus $\widetilde{\mu}$ of inclusions goes to infinity. We would like to prove that the scattered field from metascreen will be significantly amplified at some specific frequencies $\omega$ very close to 0, which are called subwavelength resonant frequencies. It is obvious that when $\omega$ is very close to the characteristic value of $\mathcal{A}$ (See Appendix \ref{sec:appa}), $\mathcal{A}$ will become ill-posed which leads to the blow-up of $\bm{\phi}$ and $\bm{\psi}$. Therefore we will begin with calculating the characteristic values of $\mathcal{A}$, prove the existence of subwavelength resonances, and then calculate the expression of $\bm{\phi}, \bm{\psi}$ followed by approximated expressions of fields in $\Omega$ and away from metascreen.

\subsection{Asymptotic analysis of Green's function and layer potential operators}\label{sec:asymptoticgreen}

Due to the complexity of problem \eqref{eq:int_repre}, it is impossible to obtain the explicit expression of $\bm{\phi}$ and $\bm{\psi}$ for an arbitrary inclusion $D$. Here we expand $\widetilde{\bm{G}}_+^{\alpha,\omega}$ asymptotically with respect to $\widetilde{\mu}$, use the leading order terms of Green's function to approximate the layer potential operators, and solve the original problem with approximated operators.

Before the calculation, we should pay additional attention to the relationship between $\widetilde{\mu}$ and quasi-periodicity $\alpha$. According to its definition, $\alpha$ is the $x_1$ component of shear wave vector $\bm{k}_s$, which implies $\alpha\leq k_s = \omega\sqrt{\rho/\widetilde{\mu}}$, hence in the subwavelength scenario $\alpha$ must converge to 0 with $\omega$. We assume here that $\alpha = \delta\omega\widetilde{\mu}^{-1/2}$ with $0<\delta<\rho$ a constant. 
We define a few variables which will be used in equations followed.
\begin{align*}
C_{\lambda,\widetilde{\mu}} = \frac{\widetilde{\mu}}{\lambda+2\widetilde{\mu}}, && \alpha_s = \sqrt{\rho-\delta^2}, && \alpha_p = \sqrt{\rho C_{\lambda,\widetilde{\mu}} - \delta^2}, \\
b_{l,s} = -\frac{\rho}{2l}, && b_{l,p} = -\frac{\rho C_{\lambda,\widetilde{\mu}}}{2l}.
\end{align*}
By direct calculation, we could verify that $\widetilde{\bm{G}}_+^{\alpha,\omega}$ could be expanded in the form of
\begin{equation}\label{eq:muinfiexpansion}
    \widetilde{\bm{G}}_+^{\alpha,\omega}(\bm{x},\bm{y}) =: \sum\limits_{k=1}^{+\infty}\frac{1}{\omega^2}\frac{\omega^{2k}}{\widetilde{\mu}^k}\widetilde{\bm{G}}^{\alpha}_{+,k}(\bm{x},\bm{y}) =:  \sum\limits_{k=1}^{+\infty}\frac{1}{\omega^2}\frac{\omega^{2k}}{\widetilde{\mu}^k}\sum\limits_{l\in 2\pi\mathbb{Z}}\widetilde{\bm{G}}_{+,k,l}^{\alpha}(\bm{x},\bm{y}),
\end{equation}
where
\begin{equation}\label{eq:g0}
    \begin{aligned}
    &\widetilde{\bm{G}}^{\alpha}_{+,1,0}(\bm{x},\bm{y}) \\
    &= \frac{1}{4\pi\rho}\begin{pmatrix}
        \rho(|x_2-y_2|-|x_2+y_2|) & -2\delta y_2(\alpha_p-\alpha_s) \\
        -2\delta y_2(\alpha_p-\alpha_s) & \rho C_{\lambda,\mu}(|x_2-y_2|-|x_2+y_2|)
    \end{pmatrix}e^{\ii\alpha (x_1-y_1)},
    \end{aligned}
\end{equation}
and for $l\neq 0$,
\begin{equation*}
    \begin{aligned}
    &\widetilde{\bm{G}}^{\alpha}_{+,1,l}(\bm{x},\bm{y})  \\
    &= \frac{e^{-l|x_2-y_2|}}{4\pi\rho}e^{\ii(\alpha+l)(x_1-y_1)}\\&\begin{pmatrix}
    b_{l,s}(1-l|x_2-y_2|) -\frac{2\delta^2}{l} + b_{l,p}(1+l|x_2-y_2|) & \ii l(b_{l,p}-b_{l,s})(x_2-y_2) \\
    \ii l(b_{l,p}-b_{l,s})(x_2-y_2) & b_{l,p}(1-l|x_2-y_2|) + b_{l,s}(1+l|x_2-y_2|) - \frac{2\delta^2}{l}
    \end{pmatrix} \\&- \frac{e^{-l|x_2+y_2|}}{4\pi\rho}e^{\ii(\alpha+l)(x_1-y_1)}\\&\begin{pmatrix}
    b_{l,s}(1-l|x_2+y_2|) -\frac{2\delta^2}{l} + b_{l,p}(1+l|x_2+y_2|) & \ii l(b_{l,p}-b_{l,s})(x_2+y_2) \\
    \ii l(b_{l,p}-b_{l,s})(x_2+y_2) & b_{l,p}(1-l|x_2+y_2|) + b_{l,s}(1+l|x_2+y_2|) - \frac{2\delta^2}{l}
    \end{pmatrix}.
    \end{aligned}
\end{equation*}
When restricting points $x,y$ inside one periodical cell, the dirac comb summation on the right-hand side of \eqref{eq:quasi-p-Green} degenerates to $\delta(x_1-y_1)\delta(x_2-y_2)\bm{I}$. Substituting \eqref{eq:muinfiexpansion} into \eqref{eq:quasi-p-Green} we get the following expansion in one periodic cell $S$ as follows
\begin{equation}\label{eq:recursiverelation}
\begin{aligned}
    &(\Delta + \nabla\nabla\cdot)\widetilde{\bm{G}}_{+,1}^{\alpha} = \delta(x_2-y_2)\delta(x_1-y_1)\bm{I}, \\
    &\omega^2(\Delta + \nabla\nabla\cdot)\widetilde{\bm{G}}_{+,k+1}^{\alpha} + \lambda\nabla\nabla\cdot\widetilde{\bm{G}}_{+,k}^{\alpha} + \rho\omega^2\widetilde{\bm{G}}_{+,k}^{\alpha} = 0,\quad k>1
\end{aligned}
\end{equation}
\begin{remark}\label{rmk:1}
We should notice that the exact expansion of $\widetilde{\bm{G}}_+^{\alpha,\omega}$ with respect to $\widetilde{\mu}$ should be
$$
    \widetilde{\bm{G}}_+^{\alpha,\omega} = \sum\limits_{p=1}^{+\infty} \widetilde{\mu}^{-p/2}\widetilde{\bm{G}}'_p + \widetilde{\mu}^{-(p+1)/2}\widetilde{\bm{G}}_p.
$$
Combining this expansion with \eqref{eq:quasi-p-Green} we have the following recursive relations
\begin{equation*}
\begin{aligned}
    &\widetilde{\bm{G}}'_1=0,\\
    &(\Delta + \nabla\nabla\cdot)\widetilde{\bm{G}}_1 = \delta(x_2-y_2)\sum\limits_{n\in\mathbb{Z}}\delta(x_1-y_1-n)\bm{I}, \\
    &(\Delta + \nabla\nabla\cdot)\widetilde{\bm{G}}_{k+2} + \lambda\nabla\nabla\cdot\widetilde{\bm{G}}_k + \rho\omega^2\widetilde{\bm{G}}_k = \delta(x_2-y_2)\sum\limits_{n\in\mathbb{Z}}\delta(x_1-y_1-n)e(k)\bm{I},\quad k>0,
\end{aligned}
\end{equation*}
where $e(k):= (\ii n\delta\omega\mu^{-1/2})^k/k!$ is the $k$-th term in expansion of $e^{\ii n\alpha}$. The last recursive equation also holds for $\widetilde{\bm{G}}'_k, k>0$. The reason for the disappearance of $\mu^{-p/2}$ terms in \eqref{eq:muinfiexpansion} is that each $\widetilde{\bm{G}}_{+,k}^{\alpha}$ there still depends on $\widetilde{\mu}$ because of the existence of $C_{\lambda,\widetilde{\mu}}$ in $\gamma_{l,p}, l\in\mathbb{N}$. The $\widetilde{\mu}^{-p/2}$ terms merge into different $\widetilde{\bm{G}}_{+,k}^{\alpha}$ in this way. It is obvious that $C_{\lambda,\widetilde{\mu}}$ is bounded when $\widetilde{\mu}\rightarrow +\infty$, and the expansion \eqref{eq:muinfiexpansion} benefits our discussion by clarifying the $\omega$ and $\widetilde{\mu}$ terms simultaneously.
\end{remark}
\begin{remark}
When calculating the recursive relation \eqref{eq:recursiverelation} we specifically assume $x$ and $y$ both lies inside one periodic cell, since this expansion will only be applied to the calculation of integrations on $\partial D$ and in $D$. It is worth mentioning that when $x,y\in\mathbb{R}^2$ are two arbitrary points, we have to use the expansion in Remark \ref{rmk:1} instead, since the $e^{\ii n\alpha}$ term in the dirac comb summation will contribute to $\widetilde{\mu}^{-p/2}$ term when we assume $\alpha = \delta\omega\widetilde{\mu}^{-1/2}$.
\end{remark}

We use $\widetilde{\bm{G}}_+^{\alpha,\omega}$(resp. $\bm{G}_+^{\alpha,\omega}$) to define $\widetilde{\mathcal{S}}^{\alpha,\omega}_+$ and $(\widetilde{\mathcal{K}}^{\alpha,\omega}_+)^{\ast}$(resp. $\mathcal{S}^{\alpha,\omega}_+$ and $(\mathcal{K}^{\alpha,\omega}_+)^{\ast}$). Notice that in $\mathcal{A}$, the left column $\widetilde{\mathcal{S}}^{\alpha,\omega}_+$ and $-\frac{1}{2}\bm{I}+(\widetilde{\mathcal{K}}_+^{\alpha,\omega})^*$ depends on asymptotic variables $\widetilde{\mu}$ and $\omega$, while the right column $-\mathcal{S}^{\alpha,\omega}_+$ and $-(\frac{1}{2}\bm{I}+(\mathcal{K}_+^{\alpha,\omega})^*)$ only depends on $\omega$. Following the decomposition of $\bm{G}_+^{\alpha,\omega}$ we have
\begin{align*}
    \widetilde{\mathcal{S}}^{\alpha,\omega}_+ = \sum\limits_{p=1}^{+\infty} \frac{\omega^{2p-2}}{\widetilde{\mu}^p}\widetilde{\mathcal{S}}^{\alpha}_{+,p}, \quad \widetilde{\mathcal{K}}^{\alpha,\omega}_+ = \sum\limits_{p=1}^{+\infty} \frac{\omega^{2p-2}}{\widetilde{\mu}^p}\widetilde{\mathcal{K}}^{\alpha}_{+,p}.
\end{align*}
Since $\widetilde{\mu}\widetilde{\bm{G}}_+^{\alpha,0} = \widetilde{\bm{G}}_{+,1}^{\alpha}$, we also have $\widetilde{\mathcal{S}}_{+,1}^{\alpha} = \widetilde{\mu}\widetilde{\mathcal{S}}_+^{\alpha,0}$ and $\widetilde{\mathcal{K}}_{+,1}^{\alpha} = \widetilde{\mu}\widetilde{\mathcal{K}}_+^{\alpha,0}$. Similar results also hold for $\bm{G}_+^{\alpha,\omega}, \mathcal{S}_+^{\alpha,\omega}$ and $\mathcal{K}_+^{\alpha,\omega}$. Furthermore, it is important to notice that $\bm{G}_{+,1}^\alpha$ is independent of choice of $\lambda,\mu$ and $\rho$, which implies that $\bm{G}_{+,1}^\alpha = \widetilde{\bm{G}}_{+,1}^\alpha$. We briefly summarize some additional useful results in the following lemma.
\begin{lemma}
    $\bm{G}_{+,1}^\alpha$ as well as $\mathcal{S}_{+,1}^\alpha$, $\mathcal{K}_{+,1}^\alpha$ is independent of $\lambda,\mu,\rho$. When regarded as operators from $L^2(\partial D)$ to $L^2(\partial D)$, the operator norm $\|\mathcal{S}_{+,1}^\alpha\|$ is of order $\mathcal{O}(1)$, while $\|\mathcal{S}_+^{\alpha,0}\| = \mathcal{O}(\mu^{-1})$, $\|\widetilde{\mathcal{S}}_+^{\alpha,0}\| = \mathcal{O}(\widetilde{\mu}^{-1})$, $\|\widetilde{\mathcal{K}}_+^{\alpha,0}\| = \mathcal{O}(1)$, $\|\mathcal{K}_+^{\alpha,0}\| = \mathcal{O}(1)$.
\end{lemma}
We can then decompose $\mathcal{A}$ as
\begin{equation}\label{eq:decomposeA}
\mathcal{A} = \mathcal{A}_0+\sum\limits_{p=1}^{+\infty}\omega^{2p-2} (\frac{1}{\widetilde{\mu}^{p}}\mathcal{A}_{p,in} + \frac{\omega^2}{\mu^{p+1}}\mathcal{A}_{p,ex}) =: \mathcal{A}_0+\sum\limits_{p=1}^{+\infty}\omega^{2p-2}\mathcal{A}_p =: \mathcal{A}_0 + \mathcal{B},
\end{equation}
where $\mathcal{A}_{p,in}$, $\mathcal{A}_{p,ex}$ is the expansion of left and right column of $\mathcal{A}$, correspondingly, in the following expressions
\begin{equation}\label{eq:adecompose}
\begin{aligned}
    \mathcal{A}_0 &= \begin{pmatrix}
    0 & -\mathcal{S}^{\alpha,0}_+ \\ -\frac{1}{2}\bm{I} + (\widetilde{\mathcal{K}}^{\alpha,0}_+)^* & -(\frac{1}{2}\bm{I}+(\mathcal{K}^{\alpha,0}_+)^*)
    \end{pmatrix}, \\
    \mathcal{A}_{p,in} &= \begin{pmatrix}
    \widetilde{\mathcal{S}}^\alpha_{+,p} & 0 \\ \frac{\omega^2}{\widetilde{\mu}}(\widetilde{\mathcal{K}}^\alpha_{+,p+1})^* & 0
    \end{pmatrix}, \qquad
    \mathcal{A}_{p,ex} = \begin{pmatrix}
    0 & -\mathcal{S}^{\alpha}_{+,p+1} \\ 0 & -(\mathcal{K}^\alpha_{+,p+1})^*
    \end{pmatrix}, \quad \mathrm{for~} p>0.
\end{aligned}
\end{equation}

  According to the discussion below Lemma 2.127 in \cite{HabibElastic}, $\mathcal{S}_+^{\alpha,0}$ might not be invertible in 2-dimensional space. As a result we impose the following assumption in this paper.
 \begin{assumption}\label{asump:invertibility}
 The shape of $D$ is chosen so that $\mathcal{S}_+^{\alpha,0}$ is invertible.
 \end{assumption}
  Although $\widetilde{\mathcal{K}}_+^{\alpha,0}$ and $\mathcal{S}_+^{\alpha,0}$ are defined with $\alpha$-quasi-periodic Green's function with Dirichlet boundary condition, we still can follow similar arguments in Section 2.15.8 of \cite{HabibElastic} to illustrate properties of interest of $\widetilde{\mathcal{K}}^{\alpha,0}_+$ and its adjoint. Define
  $$
  \Psi :=\{\bm{p}+\bm{Q}\bm{x}, \bm{p}\in\mathbb{R}^2, \bm{Q} \mathrm{~is~antisymmetric~} 2\times 2 \mathrm{~matrix}\}
  $$
  as the set of rigid motions in $\mathbb{R}^2$. We can prove in a similar way as Lemma 2.128 of \cite{HabibElastic} that $\Psi$ is the eigenspace of $\mathcal{K}^{\alpha,0}_+$ on $L^2(\partial D)$ corresponding to $\frac{1}{2}$. Since $\Psi$ is independent of $\lambda$ and $\mu$, it is immediate that $\Psi$ is also the eigenspace of $\widetilde{\mathcal{K}}^{\alpha,0}_+$.

  Let $(\cdot,\cdot)$ denote the inner product of $L^2(\partial D)\times L^2(\partial D)=:(L^2(\partial D))^2$, and $(\cdot,\cdot)_{\mathcal{H}}$ denote the inner product of $(L^2(\partial D))^2\times (L^2(\partial D))^2=:\mathcal{H}$.

It is observed that $\mathrm{dim}\Psi = 3$ since $\bm{p}$ and $\bm{Q}$ contributes two and one degrees of freedom, correspondingly. Define
$$
\bm{f}^{(1)} := \frac{1}{|\partial D|}\begin{pmatrix} 1 \\ 0\end{pmatrix}, \quad \bm{f}^{(2)} := \frac{1}{|\partial D|}\begin{pmatrix} 0 \\ 1\end{pmatrix}, \quad \bm{f}^{(3)} := \frac{1}{C_0}\begin{pmatrix}-x_2 \\ x_1\end{pmatrix}
$$
where $C_0:=\int_{\partial D} x_1^2+x_2^2 \mathrm{d}\sigma(x)$. $\bm{f}^{(i)}\in (L^2(\partial D))^2$, $i=1,2,3$ then satisfy $(\bm{f}^{(i)},\bm{f}^{(j)})=\delta_{ij}$ and therefore form an orthogonal basis of $\Psi$.
Define
$$
    H_\Psi(\partial D):=\{\bm{\varphi}\in (L^2(\partial D))^2:(\bm{\psi},\bm{\varphi}) = 0, \forall\bm{\psi}\in\Psi\}.
$$
 From Proposition \ref{prop:layerpotentialproperty} and the invertibility of $\frac{1}{2}\bm{I}-(\widetilde{\mathcal{K}}_+^{\alpha,0})^*$ on $H_\Psi(\partial D)$ from Lemma 2.127 of \cite{HabibElastic}, there exists a unique $\widetilde{\bm{\psi}}^{(i)}\in H_\Psi(\partial D)$ such that
$$
    (\frac{1}{2}\bm{I}-(\widetilde{\mathcal{K}}^{\alpha,0}_+)^*)[\widetilde{\bm{\psi}}^{(i)}] = \frac{\partial}{\partial\widetilde{\bm{\nu}}}\widetilde{S}_+^{\alpha,0}[\bm{f}^{(i)}]|_- = (-\frac{1}{2}\bm{I} + (\widetilde{\mathcal{K}}_+^{\alpha,0})^*)[\bm{f}^{(i)}].
$$
We define $\bm{\psi}^{(i)} := \widetilde{\bm{\psi}}^{(i)} + \bm{f}^{(i)}$, then it follows that $(\widetilde{\mathcal{K}}_+^{\alpha,0})^*[\bm{\psi}^{(i)}] = \frac{1}{2}\bm{\psi}^{(i)}$ and
$$
    (\bm{f}^{(i)}, \bm{\psi}^{(j)})= (\bm{f}^{(i)}, \widetilde{\bm{\psi}}^{(j)}) + (\bm{f}^{(i)}, \bm{f}^{(j)}) = \delta_{ij}.
$$

  Define $W:=\mathrm{span}\{\bm{\psi}^{(i)}\in (L^2(\partial D))^2, i=1,2,3\}$. Let
  $$
    H_W :=\{ \bm{f}\in (L^2(\partial D))^2: (\bm{f},\bm{\varphi})=0~\mathrm{for}~\mathrm{all}~\bm{\varphi}\in W\}.
  $$
  We have the following decomposition from Lemma 2.129 of \cite{HabibPhoto}.
  \begin{lemma}\label{lem:decompose}
    The following results hold.
    \begin{enumerate}
        \item Each $\bm{\varphi}\in (L^2(\partial D))^2$ is uniquely decomposed as
            $$
                \bm{\varphi} = \bm{\varphi}'+\bm{\varphi}'':=\bm{\varphi}'+\sum\limits_{j=1}^3( \bm{f}^{(j)},\bm{\varphi})\bm{\psi}^{(j)},
            $$
            and $\bm{\varphi}'\in H_\Psi$.
        \item Each $\bm{f}\in (L^2(\partial D))^2$ is uniquely decomposed as
            $$
                \bm{f} = \bm{f}' + \bm{f}'' := \bm{f}' + \sum\limits_{j=1}^3(\bm{f},\bm{\psi}^{(j)})\bm{f}^{(j)},
            $$
            and $\bm{f}'\in H_W$.
        \item $\widetilde{\mathcal{S}}_+^{\alpha,0}$ maps $W$ into $\Psi$, and $H_\Psi$ into $H_W$.
        \item $W$ is the eigenspace of $(\widetilde{\mathcal{K}}_+^{\alpha,0})^*$ corresponding to the eigenvalue $1/2$.
    \end{enumerate}
  \end{lemma}

  Due to previous discussions, we know that $\mathcal{A}_0$ is not invertible. In fact, we have
  \begin{align*}
    \mathrm{Ker}\mathcal{A}_0 &= \{(\bm{f}_1,\bm{f}_2)^\top: \bm{f}_1\in\mathrm{Ker}(-\frac{1}{2}\bm{I}+(\widetilde{\mathcal{K}}^{\alpha,0}_+)^*), \bm{f}_2=0\} \\
    &= \mathrm{span}\{\Psi_i:=(\bm{\psi}^{(i)},0)^\top, i=1,2,3\}. \\
    \mathrm{Ker}\mathcal{A}_0^* &= \{(\bm{g}_1,\bm{g}_2)^\top: \bm{g}_2\in\mathrm{Ker}(-\frac{1}{2}\bm{I}+\widetilde{\mathcal{K}}^{\alpha,0}_+), \bm{g}_1 = (\mathcal{S}_+^{\alpha,0})^{*,-1}[\bm{g}_2]\} \\
    &= \mathrm{span}\{((\mathcal{S}_+^{\alpha,0})^{*,-1}[\bm{f}^{(i)}], \bm{f}^{(i)})^\top, i=1,2,3\}.
  \end{align*}
  Define $\Psi_i^*:=((\mathcal{S}_+^{\alpha,0})^{*,-1}[\bm{f}^{(i)}], \bm{f}^{(i)})^\top/C^*_i, i=1,2,3$ where $C^*_i$ are the constants such that $(\Psi_i^*, \Psi_j^*)_{\mathcal{H}} = \delta_{ij}$.

  Therefore, similar to Lemma \ref{lem:decompose}, we can orthogonally decompose $\mathcal{H}$ using $\mathrm{Ker}\mathcal{A}_0$ as well as $\mathrm{Ker}\mathcal{A}_0^*$. We define operator $\mathcal{P}:\mathcal{H}\rightarrow \mathcal{H}$ as
  $$
  \mathcal{P}\bigg[\begin{pmatrix} \bm{f}_1 \\ \bm{f}_2 \end{pmatrix}\bigg] = \sum\limits_{i=1}^3\bigg(\begin{pmatrix} \bm{f}_1 \\ \bm{f}_2 \end{pmatrix}, \begin{pmatrix} \bm{\psi}^{(i)} \\ 0 \end{pmatrix}\bigg)_{\mathcal{H}}\Psi_i^*.
  $$
  It is immediate that $\mathcal{P}[\bm{f}] = 0$ if and only if $\bm{f}\perp \mathrm{Ker}\mathcal{A}_0$. Therefore we can construct a modified version of $\mathcal{A}_0$ which is invertible.

  \begin{lemma}\label{lem:ainvertibility}
    $\widetilde{\mathcal{A}}_0 := \mathcal{A}_0 + \mathcal{P}: \mathcal{H}\rightarrow \mathcal{H}$ is invertible. Furthermore, $\widetilde{\mathcal{A}}_0[\Psi_i]=\Psi_i^*$, $\widetilde{\mathcal{A}}_0^*[\Psi_i^*]=\Psi_i, i=1,2,3$.
  \end{lemma}

\subsection{Calculation of resonant frequencies}\label{sec:resonantfrequency}

  When regarding $\mathcal{A}$ as an operator-valued function with independent variable $\omega$, the resonant frequencies are exactly the characteristic values (defined in Appendix \ref{sec:appa}) of $\mathcal{A}$. Unfortunately, the direct calculation of characteristic values of $\mathcal{A}$ is too complicated. \cite{HabibARMA} introduced a complete method to deal with such problem, which could be summarized as follows:
  \begin{itemize}
  \item Do asymptotic expansion of $\mathcal{A}(\omega)$ to find out leading order term $\mathcal{A}_0(\omega)$;
  \item Using Rouch\'{e}'s theorem (Theorem \ref{thm:rouche}), prove that for any characteristic value $\omega$ of $\mathcal{A}$, there exists a small neighborhood $V$ of $\omega$ such that $\mathcal{A}_0$ has characteristic value $\omega_0$ in $V$ with the same multiplicity as $\omega$ (Lemma 3.8 in \cite{HabibARMA});
  \item Calculate the difference between $\omega_0$ and $\omega$ with results in Theorem 3.9 in \cite{HabibARMA}.
  \end{itemize}
  We will apply this method here. The asymptotic expansion has been done in Section \ref{sec:asymptoticgreen}. Since the steps in \cite{HabibARMA} could be directly applied here to prove the relationship between characteristic values of $\mathcal{A}$ and $\mathcal{A}_0$ without much modification, we will not repeat them for brevity. Our goal of this section is to calculate the an approximation to the characteristic value $\omega$ of $\mathcal{A}$. We will not demonstrate the differences between $\omega$ and $\omega_0$, which could still be done by following the steps in Theorem 3.9 in \cite{HabibARMA}.

  Now we suppose for the characteristic value $\omega$, $\Phi\in \mathcal{H}$ is the non-trivial solution of equation $\mathcal{A}[\Phi] = 0$. We then decompose $\Phi$ as $\Phi = \Phi_0 + \Phi_\perp$, where
  $$
  \Phi_0:=\sum_{i=1}^3 \tau_i(\bm{\psi}^{(i)},0)^\top\in\mathrm{Ker}\mathcal{A}_0, \quad \Phi_\perp \perp \mathrm{Ker}\mathcal{A}_0.
  $$
  Here $\tau_i (i=1,2,3)$ are constants to be determined. An immediate result is that $P[\Phi_\perp]=0$.

Lemma \ref{lem:ainvertibility} ensured the invertibility of $\widetilde{\mathcal{A}}_0$, and due to the fact that $\widetilde{\mu}\rightarrow\infty$ and $\omega$ is very small, Rouch\'{e}'s theorem (Chapter 1 in \cite{HabibLayer}) guarantees the invertiblity of $\widetilde{\mathcal{A}}_0+\mathcal{B}$. We have the following asymptotic expansion
\begin{equation}\label{eq:a0pbexpansion}
    (\widetilde{\mathcal{A}}_0+\mathcal{B})^{-1} = (I+\sum\limits_{p=1}^{+\infty}(-1)^p(\widetilde{\mathcal{A}}_0^{-1}\mathcal{B})^p)\widetilde{\mathcal{A}}_0^{-1}
\end{equation}
Applying it to equation $\mathcal{A}[\Phi]=0$ leads to
\begin{equation}\label{eq:omegafreq1}
    (\widetilde{\mathcal{A}}_0+\mathcal{B})^{-1}(\widetilde{\mathcal{A}}_0+\mathcal{B}-\mathcal{P})[\Phi] = 0.
\end{equation}
Since $\mathcal{P}[\Phi] = \mathcal{P}[\Phi_0] + \mathcal{P}[\Phi_\perp] = \sum_{i=1}^3\tau_i\Psi_i^*$, \eqref{eq:omegafreq1} is equivalent to
$$
    \Phi_0 + \Phi_\perp - (\widetilde{\mathcal{A}}_0+\mathcal{B})^{-1}[\sum_{i=1}^3\tau_i\Psi_i^*] = 0
$$
Taking inner product with $\Phi_0$ on both sides of \eqref{eq:omegafreq1} gets rid of $\Phi_\perp$ term and transforms the equation above to
$$
    \left(\sum\limits_{i=1}^3\tau_i\Psi_i,\sum\limits_{j=1}^3\tau_j\Psi_j\right)_\mathcal{H} - \left((I+\sum\limits_{p=1}^{+\infty}(-1)^p(\widetilde{\mathcal{A}}_0^{-1}\mathcal{B})^p)\widetilde{\mathcal{A}}_0^{-1}[\sum\limits_{i=1}^3\tau_i\Psi_i^*], \sum\limits_{j=1}^3\tau_j\Psi_j\right)_\mathcal{H} = 0.
$$
Lemma \ref{lem:ainvertibility} tells us that $\widetilde{\mathcal{A}}_0^{-1}[\Psi_i^*] = \Psi_i$, which implies that
\begin{equation}\label{eq:cal1}
    \bigg((\widetilde{\mathcal{A}}_0^{-1}\mathcal{B}-(\widetilde{\mathcal{A}}_0^{-1}\mathcal{B})^2 + \cdots)\sum\limits_{i=1}^3\tau_i\Psi_i, \sum\limits_{j=1}^3\tau_j\Psi_j\bigg)_\mathcal{H} = 0.
\end{equation}
From \eqref{eq:adecompose}, we could see that $\mathcal{A}_{k,ex}\Psi_i = 0$ holds for any $k$, therefore we only need to focus on the contribution of $\mathcal{A}_{k,in}$ when calculating \eqref{eq:cal1}. From definition of $\mathcal{B}$ (Equation \eqref{eq:decomposeA}), $((\widetilde{\mathcal{A}}_0^{-1}\mathcal{B})^k\Psi_i, \Psi_j)_\mathcal{H} = O(\widetilde{\mu}^{-2})$ holds for $k\geq 2$, therefore we can neglect all the other higher-order terms. For the leading order term of $\widetilde{\mathcal{A}}_0^{-1}\mathcal{B}$, recalling that $\mathcal{A}_1 = \frac{1}{\widetilde{\mu}}\mathcal{A}_{1,in} + \frac{\omega^2}{\mu^2}\mathcal{A}_{1,ex}$ and $(\widetilde{\mathcal{A}}_0^{-1})^*[\Psi_j] = \Psi_j^*$, we then have
\begin{equation}\label{eq:a0invb1}
    (\widetilde{\mathcal{A}}_0^{-1}\mathcal{B}\Psi_i, \Psi_j)_\mathcal{H} = \frac{1}{\widetilde{\mu}}(\widetilde{\mathcal{A}}_0^{-1}\mathcal{A}_{1,in}\Psi_i, \Psi_j)_\mathcal{H}+\mathcal{O}(\frac{1}{\widetilde{\mu}^2}) = \frac{1}{\widetilde{\mu}}(\mathcal{A}_{1,in}\Psi_i, \Psi_j^*)_\mathcal{H}+\mathcal{O}(\frac{1}{\widetilde{\mu}^2}).
\end{equation}
Since
$$
    \mathcal{A}_{1,in}\Psi_i = \begin{pmatrix} \widetilde{\mathcal{S}}_{+,1}^\alpha[\bm{\psi}^{(i)}] \\ \frac{\omega^2}{\widetilde{\mu}}(\widetilde{\mathcal{K}}_{+,2}^\alpha)^*[\bm{\psi}^{(i)}]\end{pmatrix},
$$
the leading order term in the \eqref{eq:a0invb1} takes the following form
\begin{equation}\label{eq:a1inleadingorder}
\begin{aligned}
 \frac{1}{\widetilde{\mu}}(\mathcal{A}_{1,in}\Psi_i, \Psi_j^*)_\mathcal{H} = \frac{1}{\widetilde{\mu}C^*_j}\bigg\{\left(\widetilde{\mathcal{S}}_{+,1}^\alpha[\bm{\psi}^{(i)}], (\mathcal{S}_+^{\alpha,0})^{*,-1}[\bm{f}^{(j)}]\right) + \frac{\omega^2}{\widetilde{\mu}}\left((\widetilde{\mathcal{K}}^\alpha_{+,2})^*[\bm{\psi}^{(i)}], \bm{f}^{(j)}\right)\bigg\}.
\end{aligned}
\end{equation}
Since $\widetilde{\mathcal{S}}_{+,1}^\alpha=\mathcal{S}_{+,1}^\alpha = \mu\mathcal{S}_+^{\alpha,0}$, the first term on the right-hand side of \eqref{eq:a1inleadingorder} is
$$
    \left(\widetilde{\mathcal{S}}_{+,1}^\alpha[\bm{\psi}^{(i)}], (\mathcal{S}_+^{\alpha,0})^{*,-1}[\bm{f}^{(j)}]\right) = (\mu\bm{\psi}^{(i)}, \bm{f}^{(j)}) = \mu\delta_{ij}.
$$
Now we calculate  $\frac{\omega^2}{\widetilde{\mu}}((\widetilde{\mathcal{K}}^\alpha_{+,2})^*[\bm{\psi}^{(i)}], \bm{f}^{(j)})$. From recursive relations \eqref{eq:recursiverelation} we know that
$$
    \frac{1}{\widetilde{\mu}}\mathcal{L}^{\widetilde{\lambda},\widetilde{\mu}}\widetilde{\bm{G}}_{+,2}^\alpha = (\Delta + \nabla\nabla\cdot + \frac{\widetilde{\lambda}}{\widetilde{\mu}}\nabla\nabla\cdot)\widetilde{\bm{G}}_{+,2}^\alpha = -(\frac{\widetilde{\lambda}}{\omega^2}\nabla\nabla\cdot + \widetilde{\rho})\widetilde{\bm{G}}_{+,1}^\alpha + O(\frac{1}{\widetilde{\mu}}).
$$
Notice that $\bm{f}^{(j)}$ satisfies
\begin{equation*}
\left\{
\begin{aligned}
&\mathcal{L}^{\lambda,\mu}\bm{f}^{(j)}=0 \quad \mathrm{in~} D, \\
&\frac{\partial \bm{f}^{(j)}}{\partial \bm{\nu}}=0 \quad \mathrm{on~} \partial D,
\end{aligned}
\right.
\end{equation*}
for any $(\lambda,\mu)$ pair. We then have
\begin{align*}
    \frac{1}{\widetilde{\mu}}\widetilde{\mathcal{K}}_{+,2}^\alpha[\bm{f}^{(j)}] &= \frac{1}{\widetilde{\mu}}\int_{\partial D}\frac{\partial \widetilde{\bm{G}}_{+,2}^\alpha(x,y)}{\partial\bm{\widetilde{\nu}}_y}\bm{f}^{(j)}(y)\mathrm{d}\sigma(y) = \frac{1}{\widetilde{\mu}}\int_D \mathcal{L}^{\widetilde{\lambda},\widetilde{\mu}}\widetilde{\bm{G}}_{+,2}^\alpha(x,y)\bm{f}^{(j)}(y)\mathrm{d}y\\
    & = -\int_D \left(\frac{\widetilde{\lambda}}{\omega^2}\nabla_y\nabla_y\cdot\widetilde{\bm{G}}_{+,1}^\alpha(x,y)+\widetilde{\rho}\widetilde{\bm{G}}_{+,1}^\alpha(x,y)\right)\bm{f}^{(j)}(y)\mathrm{d}y + O(\frac{1}{\widetilde{\mu}}) \\
    &= \int_D \frac{\widetilde{\lambda}}{\omega^2}\Delta_y\widetilde{\bm{G}}_{+,1}^\alpha(x,y) \bm{f}^{(j)}(y)\mathrm{d}y -\widetilde{\rho}\int_D\widetilde{\bm{G}}_{+,1}^\alpha(x,y)\bm{f}^{(j)}(y)\mathrm{d}y + O(\frac{1}{\widetilde{\mu}}).
\end{align*}
Notice that $\widetilde{\bm{G}}_{+,1}^\alpha$ satisfies homogeneous Lam\'{e} system with $\lambda = 0, \mu = 1,\omega=0$, which means the corresponding boundary contraction is
$$
\frac{\partial\widetilde{\bm{G}}_{+,1}^\alpha}{\partial\bm{\nu}_{1,0,y}} = (\nabla_y\widetilde{\bm{G}}_{+,1}^\alpha + \nabla_y(\widetilde{\bm{G}}_{+,1}^\alpha)^\top)\bm{\nu}_y = \partial_{\bm{\nu}_y}\widetilde{\bm{G}}_{+,1}^\alpha + \partial_{\bm{\nu}_y}(\widetilde{\bm{G}}_{+,1}^\alpha)^\top = 2\partial_{\bm{\nu}_y}\widetilde{\bm{G}}_{+,1}^\alpha.
$$
Since $x\in\partial D$, calculation then follows
\begin{align*}
    \int_D \Delta_y \widetilde{\bm{G}}_{+,1}^\alpha(x,y)\bm{f}^{(j)}(y)\mathrm{d}y &= \int_{\partial D}\partial_{\bm{\nu}_y}\widetilde{\bm{G}}_{+,1}^\alpha(x,y)\bm{f}^{(j)}(y)\mathrm{d}\sigma(y) = \frac{1}{2}\int_{\partial D}\frac{\partial \widetilde{\bm{G}}_{+,1}^\alpha(x,y)}{\partial \bm{\nu}_{1,0,y}}\bm{f}^{(j)}(y)\mathrm{d}\sigma(y)\\
     &= \frac{1}{2}\int_D \mathcal{L}^{1,0}\widetilde{\bm{G}}_{+,1}^\alpha(x,y)\bm{f}^{(j)}(y)\mathrm{d}\sigma(y)=0.
\end{align*}
Therefore the second term in $(\mathcal{A}_{1,in}\Psi_i,\Psi_j^*)_\mathcal{H}$ is
\begin{align*}
    \frac{\omega^2}{\widetilde{\mu}}((\widetilde{\mathcal{K}}_{+,2}^\alpha)^*[\bm{\psi}^{(i)}],\bm{f}^{(j)}) &= \frac{\omega^2}{\widetilde{\mu}}(\bm{\psi}^{(i)}, \widetilde{\mathcal{K}}_{+,2}^\alpha[\bm{f}^{(j)}]) \\
    &= -\widetilde{\rho}\omega^2\int_{\partial D}\bm{\psi}^{(i)}(x)\int_D\widetilde{\bm{G}}_{+,1}^\alpha(x,y)\bm{f}^{(j)}(y)\mathrm{d}y\mathrm{d}\sigma(x) + O(\frac{1}{\widetilde{\mu}}) \\
    &= -\widetilde{\rho}\omega^2\int_D \widetilde{\mathcal{S}}_{+,1}^\alpha[\bm{\psi}^{(i)}](y)\bm{f}^{(j)}(y)\mathrm{d}y + O(\frac{1}{\widetilde{\mu}}).
\end{align*}
Define
\begin{equation*}
\mathcal{C} = \begin{pmatrix}\frac{1}{C_1^*} & & \\ & \frac{1}{C_2^*} & \\ & & \frac{1}{C_3^*}\end{pmatrix}, \quad \mathcal{M} = \begin{pmatrix} m_{11} & m_{12} & m_{13} \\ m_{21} & m_{22} & m_{23} \\ m_{31} & m_{32} & m_{33} \end{pmatrix}
\end{equation*}
where $m_{ij} := \int_D \widetilde{\mathcal{S}}_{+,1}^{\alpha}[\bm{\psi}^{(i)}](y)\bm{f}^{(j)}(y)\mathrm{d}y$. We can then summarize the calculation above as
\begin{equation}\label{eq:a0inva1}
    (\widetilde{\mathcal{A}}_0^{-1}\mathcal{A}_1\sum\limits_{i=1}^3\tau_i\Psi_i, \sum\limits_{i=1}^3\tau_j\Psi_j)_\mathcal{H} = \bm{\tau}^\top(\mu\bm{I}-\widetilde{\rho}\omega^2\mathcal{M}+O(\frac{1}{\widetilde{\mu}^2}))\mathcal{C}\bm{\tau}.
\end{equation}
The characteristic values of $\mathcal{A}$ can then be approximated by the eigenvalues of $\mu\bm{I} - \widetilde{\rho}\omega^2\mathcal{M}$. To conclude:
\begin{theorem}\label{thm:resonantfrequency}
    Let $m_i, i=1,2,3$ be the eigenvalue of $\mathcal{M}$. Then the resonant frequencies are
    $$
    \omega_i = \sqrt{\frac{\mu}{\widetilde{\rho} m_i}} + O(\widetilde{\mu}^{-1}).
    $$
    Since $m_i$ only depends on the shape of $D$, it is immediate that the subwavelength resonance occurs \emph{if and only if} $\mu/\widetilde{\rho}\rightarrow 0$.
\end{theorem}


\subsection{Near-field and far-field calculation}
In the scenario of sub-wavelength resonance, i.e. $\omega\rightarrow 0$ and $\widetilde{\mu}\rightarrow \infty$, we are able to carry out further calculations of approximated near field and far field based on asymptotic expansions of $\omega$ and $\widetilde{\mu}$. Assuming now the incident wave is P-wave, we want to solve
\begin{equation}\label{eq:psimain}
(\widetilde{\mathcal{A}}_0+\mathcal{B}-\mathcal{P})\Phi = \bm{F}.
\end{equation}
$\bm{F}$ is defined in \eqref{eq:Fdefinition}. Same as previous subsection, we decompose $\Phi = \Phi_0 + \Phi_\perp$, where $\Phi_0 = \sum_{i=1}^3\tau_i\Psi_i$. The first goal of this subsection is to determine $\tau_i$ and $\Phi_\perp$.

\subsubsection{Determine $\tau_i$.}

By applying $(\widetilde{\mathcal{A}}_0+\mathcal{B})^{-1}$ to both sides of \eqref{eq:psimain} we get
\begin{equation}\label{eq:psimain1}
    \Phi_0 + \Phi_\perp - (\widetilde{\mathcal{A}}_0+\mathcal{B})^{-1}\mathcal{P}\Phi = (\widetilde{\mathcal{A}}_0+\mathcal{B})^{-1}\bm{F},
\end{equation}
where $\bm{F} = (\bm{u}_P^+, \partial\bm{u}_P^+/\partial\bm{\nu})^\top$. The contraction $\partial\bm{\nu}$ is taking on the outward boundary of $\partial D$, hence is independent of $\widetilde{\mu}$. Taking inner product with $\Psi_j$ on both sides and following the same calculation steps as \eqref{eq:omegafreq1}-\eqref{eq:cal1}, we get
\begin{equation}\label{eq:psijinnerproduct}
    \bigg((\widetilde{\mathcal{A}}_0^{-1}\mathcal{B}-(\widetilde{\mathcal{A}}_0^{-1}\mathcal{B})^2 + \cdots)\sum\limits_{i=1}^3\tau_i\Psi_i, \Psi_j\bigg)_\mathcal{H} = ((\widetilde{\mathcal{A}}_0+\mathcal{B})^{-1}\bm{F}, \Psi_j)_\mathcal{H}.
\end{equation}
The right-hand side of \eqref{eq:psijinnerproduct} can be expanded as
\begin{equation}\label{eq:a0binv1}
    ((\widetilde{\mathcal{A}}_0+\mathcal{B})^{-1}\bm{F}, \Psi_j)_\mathcal{H} = (\widetilde{\mathcal{A}}_0^{-1}\bm{F}-\widetilde{\mathcal{A}}_0^{-1}\mathcal{A}_1\widetilde{\mathcal{A}}_0^{-1}\bm{F}, \Psi_j)_\mathcal{H} + \mathcal{O}(\omega^2)+\mathcal{O}(\frac{1}{\widetilde{\mu}^2}).
\end{equation}
By carefully checking the definition of $\mathcal{A}_1$, we notice that
\begin{equation}\label{eq:a1expansion}
    \mathcal{A}_1 = \frac{1}{\widetilde{\mu}}\begin{pmatrix}\widetilde{\mathcal{S}}_{+,1}^\alpha & 0 \\ 0 & 0\end{pmatrix} + \frac{\omega^2}{\widetilde{\mu}^2}\begin{pmatrix}0 & 0 \\ (\widetilde{\mathcal{K}}_{+,2}^\alpha)^* & 0\end{pmatrix} - \frac{\omega^2}{\mu^2}\begin{pmatrix} 0 & \mathcal{S}_{+,2}^\alpha \\ 0 & (\mathcal{K}_{+,2}^\alpha)^*\end{pmatrix} = \frac{1}{\widetilde{\mu}}\begin{pmatrix}\widetilde{\mathcal{S}}_{+,1}^\alpha & 0 \\ 0 & 0\end{pmatrix} + \mathcal{O}(\omega^2).
\end{equation}
Therefore, the leading order term of $(\widetilde{\mathcal{A}}_0^{-1}\mathcal{A}_1\widetilde{\mathcal{A}}_0^{-1}\bm{F},\Psi_j)$ is
$$
    \frac{1}{\widetilde{\mu}}\left(\begin{pmatrix}\widetilde{\mathcal{S}}_{+,1}^\alpha & 0 \\ 0 & 0\end{pmatrix}\widetilde{\mathcal{A}}_0^{-1}\bm{F}, \Psi_j^*\right)_\mathcal{H} = \frac{1}{\widetilde{\mu}}\left(\widetilde{\mathcal{A}}_0^{-1}\bm{F}, \begin{pmatrix}(\widetilde{\mathcal{S}}_{+,1}^\alpha)^* & 0 \\ 0 & 0\end{pmatrix}\begin{pmatrix} (\mathcal{S}_+^{\alpha,0})^{*,-1}[\bm{f}^{(j)}] \\ \bm{f}^{(j)}\end{pmatrix}\right)_\mathcal{H}= \frac{\mu}{\widetilde{\mu}}(\widetilde{\mathcal{A}}_0^{-1}\bm{F}, \Psi_j)_\mathcal{H}.
$$
The equation \eqref{eq:a0binv1} could then be reformulated as
\begin{equation}\label{eq:a0binv2}
    ((\widetilde{\mathcal{A}}_0+B)^{-1}\bm{F}, \Psi_j)_\mathcal{H} = \left(1-\frac{\mu}{\widetilde{\mu}}\right)(\widetilde{\mathcal{A}}_0^{-1}\bm{F}, \Psi_j) + \mathcal{O}(\omega^2)+\mathcal{O}(\frac{1}{\widetilde{\mu}^2}) =: \left(1-\frac{\mu}{\widetilde{\mu}}\right)\mathcal{K}_j + \mathcal{O}(\omega^2)+\mathcal{O}(\frac{1}{\widetilde{\mu}^2}).
\end{equation}
Define $\mathcal{K}:=(\mathcal{K}_1, \mathcal{K}_2, \mathcal{K}_3)^\top$. Equation \eqref{eq:a0inva1} tells us that
$$
((\widetilde{\mathcal{A}}_0+B)^{-1}[\sum\tau_i\Psi_i^*], \Psi_j) = \mathcal{C}(\mu\bm{I}-\widetilde{\rho}\omega^2\mathcal{M})\bm{\tau} + O(\frac{1}{\widetilde{\mu}^2}).
$$
Plugging \eqref{eq:a0binv2} and this equation into \eqref{eq:psijinnerproduct}, we obtain
$$
\bm{\tau} = (1-\frac{\mu}{\widetilde{\mu}})(\mu\bm{I}-\widetilde{\rho}\omega^2\mathcal{M})^{-1}\mathcal{C}^{-1}\mathcal{K}+\mathcal{O}(\omega^2)+\mathcal{O}(\widetilde{\mu}^{-2}).
$$
To achieve an estimation of $\mathcal{K}_j$, we need to calculate
$$
(\widetilde{\mathcal{A}}_0^{-1}F, \Psi_j)_\mathcal{H} = (F, \Psi_j^*)_\mathcal{H} = \left(\begin{pmatrix}\bm{u}_+^P \\ \frac{\partial \bm{u}_+^P}{\partial \bm{\nu}}\end{pmatrix}, \begin{pmatrix}(\mathcal{S}^{\alpha,0}_+)^{*,-1}[\bm{f}^{(j)}] \\ \bm{f}^{(j)}\end{pmatrix}\right)_\mathcal{H}.
$$
It then suffices to calculate $(\bm{u}_+^P, (\mathcal{S}^{\alpha,0}_+)^{*,-1}[\bm{f}^{(j)}])$ and $(\partial \bm{u}_+^P/\partial\bm{\nu}, \bm{f}^{(j)})$.
\begin{itemize}
\item \emph{$(\bm{u}_+^P, (\mathcal{S}^{\alpha,0}_+)^{*,-1}[\bm{f}^{(j)}])$ term}.

Since $\omega\rightarrow 0$, we could expand definition equation \eqref{eq:upp} of $\bm{u}_P^+$ and get that
$$
    \bm{u}_+^P(\bm{x}) = 2\ii\omega\sqrt{\frac{\rho}{\lambda+2\mu}}\theta_2x_2\bm{\theta} + \mathcal{O}(\omega^2).
$$
Hence
$$
    (\bm{u}_+^P, (\mathcal{S}_+^{\alpha,0})^{*,-1}[\bm{f}^{(j)}]) = 2\ii\omega\sqrt{\frac{\rho}{\lambda+2\mu}}\theta_2((\mathcal{S}_+^{\alpha,0})^{-1}[x_2\bm{\theta}], \bm{f}^{(j)}) + \mathcal{O}(\omega^2).
$$
To proceed, we need the following lemma.
\begin{lemma}\label{lem:htheta}
Define
$$
\bm{h}(\bm{\theta}):=\frac{\partial (x_2\bm{\theta})}{\partial\bm{\nu}} = \theta_2(\lambda\nu_1, (\lambda+2\mu)\nu_2)^\top.
$$
Then
    $$
    (\mathcal{S}_+^{\alpha,0})^{-1}[x_2\bm{\theta}] = (-\frac{1}{2}\bm{I}+(\mathcal{K}^{\alpha,0}_+)^*)^{-1}[\bm{h}(\bm{\theta})] + \sum\limits_{i=1}^3 ( x_2\bm{\theta},\bm{\psi}^{(i)})(\mathcal{S}_+^{\alpha,0})^{-1}[\bm{f}^{(i)}].
    $$
\end{lemma}
\begin{proof}
Define $\bm{f}_\theta=x_2\bm{\theta} - \sum_{i=1}^3(x_2\bm{\theta}, \bm{\psi}^{(i)}) \bm{f}^{(i)}$. The construction of $\bm{f}_\theta$ immediately implies that $\bm{f}_\theta\in H_W(\partial D)$. Recalling that $\bm{f}^{(i)}\in\Psi$ are solutions of inner Neumann problem, we can then observe that
\begin{equation*}
\left\{
\begin{aligned}
&\mathcal{L}^{\lambda,\mu}\bm{f}_\theta=0 \quad \mathrm{in}~D, \\
&\frac{\partial \bm{f}_\theta}{\partial\bm{\nu}}|_- = \bm{h}(\bm{\theta}) \quad\mathrm{on}~\partial D.
\end{aligned}
\right.
\end{equation*}
It follows that $\bm{h}(\bm{\theta})\in H_\Psi(\partial D)$, which means $\bm{h}(\bm{\theta})\in\mathrm{Im}(-\bm{I}/2+(\mathcal{K}_+^{\alpha,0})^*)$. Define $\bm{\eta}_-:=(-\bm{I}/2+(\mathcal{K}^{\alpha,0}_+)^*)^{-1}[\bm{h}(\bm{\theta})]$. Then $\bm{\eta}_-\in H_\Psi(\partial D)$ and $\mathcal{S}_+^{\alpha,0}[\bm{\eta}_-] = \bm{f}_\theta$ holds, since we have already ensured that $\bm{f}_\theta\in H_W(\partial D)$. The result then follows by applying $(\mathcal{S}_+^{\alpha,0})^{-1}$ on both sides of definition equation of $\bm{f}_\theta$.
\end{proof}
From Lemma \ref{lem:decompose} we know that $\mathcal{S}_+^{\alpha,0}$ maps $W$ into $\Psi$. Recall that $\widetilde{\mathcal{S}}_{+,1}^\alpha$ is independent of $\lambda,\mu$ and satisfies $\widetilde{\mathcal{S}}_{+,1}^\alpha = \mu\mathcal{S}_+^{\alpha,0}$. We assume
$$
\widetilde{\mathcal{S}}_{+,1}^\alpha[\bm{\psi}^{(i)}] =: \sum_{j=1}^3 s_{ij}\bm{f}^{(j)}
$$
and define $\bm{S}=(s_{ij})_{3\times 3}$. The invertibility of $\widetilde{\mathcal{S}}_{+,1}^{\alpha}$ indicates $\bm{S}$ is invertible. Denote $\bm{S}^{-1} = (s_{ij}^{-1})_{3\times 3}$, where $\sum_{j=1}^3 s_{ij}^{-1}s_{jk} = \delta_{ik}$. Then $(\widetilde{\mathcal{S}}_{+,1}^\alpha)^{-1}[\bm{f}^{(i)}] = \sum_{j=1}^3 s_{ij}^{-1}\bm{\psi}^{(j)}$. As a result,
\begin{align*}
    \left((\mathcal{S}_+^{\alpha,0})^{-1}[x_2\bm{\theta}], \bm{f}^{(i)}\right) =& \left((\mathcal{S}_+^{\alpha,0})^{-1}[\bm{f}_\theta] + \sum\limits_{j=1}^3( x_2\bm{\theta},\bm{\psi}^{(j)})(\mathcal{S}_+^{\alpha,0})^{-1}[\bm{f}^{(j)}], \bm{f}^{(i)}\right) \\
    =& \sum\limits_{j=1}^3 (x_2\bm{\theta},\bm{\psi}^{(j)})\sum\limits_{k=1}^3 (\mu s_{jk}^{-1}\bm{\psi}^{(k)}, \bm{f}^{(i)}) = \sum\limits_{j=1}^3(x_2\bm{\theta}, \bm{\psi}^{(j)})\mu s_{ji}^{-1}.
\end{align*}
The second equality holds because $(\mathcal{S}_+^{\alpha,0})^{-1}[\bm{f}_\theta]\in H_\Psi(\partial D)$. We then have
\begin{align*}
\left(\bm{u}_+^P, (\mathcal{S}_+^{\alpha,0})^{*,-1}[\bm{f}^{(j)}]\right) = 2\ii\omega\sqrt{\frac{\rho}{\lambda+2\mu}}\theta_2 \sum\limits_{i=1}^3(x_2\bm{\theta}, \bm{\psi}^{(i)})\mu s_{ij}^{-1} + \mathcal{O}(\omega^2).
\end{align*}
\item \emph{$(\partial \bm{u}_+^P/\partial\bm{\nu}, \bm{f}^{(i)})$ term.}

By expanding \eqref{eq:upnu1}, \eqref{eq:upnu21} and \eqref{eq:upnu22} with respect to $\omega$, we know that
\begin{align*}
    \left(\frac{\partial \bm{u}_+^P}{\partial \bm{\nu}}\right)_1 &= 2\lambda\omega\sqrt{\frac{\rho}{\lambda+2\mu}}\nu_1 + 2\ii\mu\omega\sqrt{\frac{\rho}{\lambda+2\mu}}\theta_1\theta_2\nu_2+ \mathcal{O}(\omega^2), \\
    \left(\frac{\partial\bm{u}_+^P}{\partial\bm{\nu}}\right) &= 2\lambda\omega\sqrt{\frac{\rho}{\lambda+2\mu}}\nu_2 + 2\ii\mu\omega\sqrt{\frac{\rho}{\lambda+2\mu}}\theta_1\theta_2\nu_1 + 4\ii\mu\omega\sqrt{\frac{\rho}{\lambda+2\mu}}\theta_2^2\nu_2+\mathcal{O}(\omega^2).
\end{align*}
So
\begin{align*}
    \left(\frac{\partial \bm{u}_+^P}{\partial \bm{\nu}}, \bm{f}^{(i)}\right)=& 2\lambda\omega\sqrt{\frac{\rho}{\lambda+2\mu}}(\bm{\nu}, \bm{f}^{(i)}) + 2\ii\mu\omega\sqrt{\frac{\rho}{\lambda+2\mu}}\theta_1\theta_2((\nu_2,\nu_1)^\top, \bm{f}^{(i)}) \\
    &+ 4\ii\mu\omega\sqrt{\frac{\rho}{\lambda+2\mu}}\theta_2^2((0,\nu_2)^\top, \bm{f}^{(i)})+ \mathcal{O}(\omega^2) .
\end{align*}
From definition of $\bm{f}^{(i)}$ we could verify that $\nabla\cdot \bm{f}^{(i)}=0$. Hence by divergence theorem we have $(\bm{\nu}, \bm{f}^{(i)}) = 0, i=1,2,3$. Similarly,
\begin{align*}
    \int_{\partial D}\begin{pmatrix}\nu_2 \\ \nu_1 \end{pmatrix}\cdot\begin{pmatrix}\bm{f}^{(i)}_1 \\ \bm{f}^{(i)}_2 \end{pmatrix}\mathrm{d}\sigma = \int_D\nabla\cdot\begin{pmatrix}\bm{f}^{(i)}_2 \\ \bm{f}^{(i)}_1\end{pmatrix}\mathrm{d}x = 0, \\
    \int_{\partial D}\begin{pmatrix}0 \\ \nu_2 \end{pmatrix}\cdot\begin{pmatrix}\bm{f}^{(i)}_1 \\ \bm{f}^{(i)}_2 \end{pmatrix}\mathrm{d}\sigma = \int_D\nabla\cdot\begin{pmatrix}0 \\ \bm{f}^{(i)}_2\end{pmatrix}\mathrm{d}x = 0.
\end{align*}
This implies $(\partial \bm{u}_+^P/\partial\bm{\nu}, \bm{f}^{(i)})=\mathcal{O}(\omega^2)$.
\end{itemize}
To conclude, the coefficient vector $\bm{\tau}$ takes the following form:
\begin{equation}\label{eq:alphaexpression}
\begin{aligned}
   \bm{\tau} =& 2\ii\omega\sqrt{\frac{\rho}{\lambda+2\mu}}\mu(1-\frac{\mu}{\widetilde{\mu}})\theta_2(\mu\bm{I}-\widetilde{\rho}\omega^2\mathcal{M})^{-1}\sum\limits_{j=1}^3(x_2\bm{\theta}, \bm{\psi}^{(j)})\begin{pmatrix}c_1s_{j1}^{-1} \\ c_2s_{j2}^{-1} \\ c_3s_{j3}^{-1}\end{pmatrix}\\&+\mathcal{O}(\omega^2)+\mathcal{O}(\widetilde{\mu}^{-2}).
\end{aligned}
\end{equation}
\begin{remark}\label{rmk:order}
Here we would like to estimate the order of $\bm{\alpha}$. From previous discussions we know that
$$
    \mathcal{C}(\mu\bm{I}-\widetilde{\rho}\omega^2\mathcal{M})\bm{\tau} = (1-\frac{\mu}{\widetilde{\mu}})\mathcal{K} + \mathcal{O}(\omega^2) + \mathcal{O}(\frac{1}{\widetilde{\mu}^2}).
$$
From definitions it is obvious that $\mathcal{C} = \mathcal{O}(1)$, $\mathcal{K} = \mathcal{O}(1)$, so the right-hand side is $\mathcal{O}(1)$. We can then conclude that $\bm{\tau} = \mathcal{O}(\mu^{-1})$. In case that $\mu = \mathcal{O}(\widetilde{\mu}^{-1})$, we then have $\bm{\tau} = \mathcal{O}(\widetilde{\mu})$.
\end{remark}

\subsubsection{Determining $\Phi_\perp$.}

Now we start the calculation of $\Phi_\perp$. Similar asymptotic analysis of \eqref{eq:psimain1} leads to
\begin{equation}\label{eq:psiperp}
\Phi_\perp = \widetilde{\mathcal{A}}_0^{-1}\bm{F} - \widetilde{\mathcal{A}}_0^{-1}\mathcal{A}_1\widetilde{\mathcal{A}}_0^{-1}\bm{F} - \sum\limits_{i=1}^3\tau_i\widetilde{\mathcal{A}}_0^{-1}\mathcal{A}_1[\Psi_i] + \mathcal{O}(\omega^2) + \mathcal{O}(\frac{1}{\widetilde{\mu}^2}).
\end{equation}
\begin{itemize}
  \item \emph{$\widetilde{\mathcal{A}}_0^{-1}\bm{F}$ term}.

  Recall that $\widetilde{\mathcal{A}}_0 = \mathcal{A}_0+\mathcal{P}$. Denote $\beta_i:= (\bm{F},\Psi_i^*)_\mathcal{H}$ and $\bm{F}_\perp:= \bm{F}-\sum_{i=1}^3\beta_i\Psi_i^*$, then $\bm{F}_\perp$ is orthogonal to $\mathrm{Ker}\mathcal{A}_0^*$, and therefore belongs to $\mathrm{Im}\mathcal{A}_0$. We then have
  $$
    \widetilde{\mathcal{A}}_0^{-1}\bm{F} = \mathcal{A}_0^{-1}\bm{F}_\perp + \sum\limits_{i=1}^3\beta_i\Psi_i.
  $$
  From the asymptotic expansion of $\bm{F}$, we can estimate $\beta_i$ by
  \begin{align*}
    \beta_i &= (\bm{F},\Psi_i^*)_\mathcal{H} = \left(\begin{pmatrix}2\ii\omega\sqrt{\frac{\rho}{\lambda+2\mu}}\theta_2x_2\bm{\theta}+\mathcal{O}(\omega^2)\\ \partial \bm{u}_+^P/\partial\bm{\nu} \end{pmatrix}, \begin{pmatrix}(\mathcal{S}_+^{\alpha,0})^{*,-1}[\bm{f}^{(i)}] \\ \bm{f}^{(i)}\end{pmatrix}\right)_\mathcal{H} \\
    &= \left(2\ii\omega\sqrt{\frac{\rho}{\lambda+2\mu}}\theta_2x_2\bm{\theta}, (\mathcal{S}_+^{\alpha,0})^{*,-1}[\bm{f}^{(i)}]\right) + \left(\frac{\partial\bm{u}_+^P}{\partial\bm{\nu}}, \bm{f}^{(i)}\right) +\mathcal{O}(\omega^2\mu) \\
    &= 2\ii\omega\sqrt{\frac{\rho}{\lambda+2\mu}}\theta_2((\mathcal{S}_+^{\alpha,0})^{-1}[x_2\bm{\theta}], \bm{f}^{(i)}) +\mathcal{O}(\omega^2) + \mathcal{O}(\omega^2\mu) \\&= 2\ii\omega\mu\theta_2\sqrt{\frac{\rho}{\lambda+2\mu}}\sum\limits_{j=1}^3(x_2\bm{\theta}, \bm{\psi}^{(j)})s_{ji}^{-1}+\mathcal{O}(\omega^2) + \mathcal{O}(\omega^2\mu).
  \end{align*}
  The second last equality holds because of results in last section. The last equality holds because of Lemma \ref{lem:htheta} and the fact that $(-\frac{1}{2}\bm{I}+(\mathcal{K}^{\alpha,0}_+)^*)^{-1}[\bm{h}(\bm{\theta})]\in H_\Psi(\partial D)$.

  Denote $\bm{F}_\perp =: (\bm{F}_{\perp,1}, \bm{F}_{\perp,2})^\top$, then
  \begin{align*}
    \bm{F}_{\perp,1} &= 2\ii\omega\sqrt{\frac{\rho}{\lambda+2\mu}}\theta_2\left(x_2\bm{\theta} - \mu\sum\limits_{i=1}^3(x_2\bm{\theta}, \bm{\psi}^{(j)})\sum\limits_{j=1}^3 s_{ji}^{-1}(\mathcal{S}_+^{\alpha,0})^{*,-1}[\bm{f}^{(i)}]\right)+\mathcal{O}(\omega^2) +\mathcal{O}(\omega^2\mu),\\
    \bm{F}_{\perp,2} &= 2\lambda\omega\sqrt{\frac{\rho}{\lambda+2\mu}}\bm{\nu} + 4\ii\mu\omega\sqrt{\frac{\rho}{\lambda+2\mu}}\theta_2^2\begin{pmatrix}0 \\ \nu_2\end{pmatrix} \\&+ 2\ii\omega\mu\sqrt{\frac{\rho}{\lambda+2\mu}}\theta_2\left(\theta_1\begin{pmatrix}\nu_2 \\ \nu_1\end{pmatrix} - \sum\limits_{i=1}^3\sum\limits_{j=1}^3(x_2\bm{\theta}, \bm{\psi}^{(j)})s_{ji}^{-1}\bm{f}^{(i)}\right) + \mathcal{O}(\omega^2) + \mathcal{O}(\omega^2\mu).
  \end{align*}
  For $\mathcal{A}_0^{-1}\bm{F}_\perp$, we estimate each item separately. We know that
  $$
    \mathcal{A}_0^{-1} = \begin{pmatrix}
    -(-\frac{1}{2}\bm{I} + (\widetilde{\mathcal{K}}_+^{\alpha,0})^*)^{-1}(\frac{1}{2}\bm{I} + (\mathcal{K}_+^{\alpha,0})^*)(\mathcal{S}_+^{\alpha,0})^{-1} & (-\frac{1}{2}\bm{I}+(\widetilde{\mathcal{K}}_+^{\alpha,0})^*)^{-1} \\ -(\mathcal{S}_+^{\alpha,0})^{-1} & 0
    \end{pmatrix}.
  $$
  So
  \begin{align*}
    (\mathcal{A}_0^{-1}\bm{F}_\perp)_1 &= -(-\frac{1}{2}\bm{I} + (\widetilde{\mathcal{K}}_+^{\alpha,0})^*)^{-1}(\frac{1}{2}\bm{I} + (\mathcal{K}_+^{\alpha,0})^*)(\mathcal{S}_+^{\alpha,0})^{-1}[\bm{F}_{\perp,1}] + (-\frac{1}{2}\bm{I}+(\widetilde{\mathcal{K}}_+^{\alpha,0})^*)^{-1}[\bm{F}_{\perp,2}], \\
    (\mathcal{A}_0^{-1}\bm{F}_\perp)_2 &= -(\mathcal{S}_+^{\alpha,0})^{-1}[\bm{F}_{\perp,1}].
  \end{align*}
  Use the definition and result in Lemma \ref{lem:htheta}, we have
  \begin{align*}
    &(\mathcal{A}_0^{-1}\bm{F}_\perp)_1 = -(-\frac{1}{2}\bm{I} + (\widetilde{\mathcal{K}}_+^{\alpha,0})^*)^{-1}(\frac{1}{2}\bm{I}+(\mathcal{K}_+^{\alpha,0})^*)(\mathcal{S}_+^{\alpha,0})^{-1}\bm{F}_{\perp,1} \\&+ 2\lambda\omega\sqrt{\frac{\rho}{\lambda+2\mu}}(-\frac{1}{2}\bm{I}+(\widetilde{\mathcal{K}}_+^{\alpha,0})^*)^{-1}[\bm{\nu}] +\mathcal{O}(\omega\mu) + \mathcal{O}(\omega^2) \\&= -2\ii\omega\sqrt{\frac{\rho}{\lambda+2\mu}}\theta_2(-\frac{1}{2}\bm{I} + (\widetilde{\mathcal{K}}_+^{\alpha,0})^*)^{-1}
    \\&\left((-\frac{1}{2}\bm{I}+(\mathcal{K}_+^{\alpha,0})^*)^{-1}[\bm{h(\theta)}] - \bm{h(\theta)} - (\frac{1}{2}\bm{I}+(\mathcal{K}_+^{\alpha,0})^*)(\mathcal{S}_+^{\alpha,0})^{-1}\sum\limits_{i=1}^3(x_2\bm{\theta},\bm{\psi}^{(i)})[\bm{f}^{(i)}]\right)\\ &+ 2\lambda\omega\sqrt{\frac{\rho}{\lambda+2\mu}}(-\frac{1}{2}\bm{I}+(\widetilde{\mathcal{K}}_+^{\alpha,0})^*)^{-1}[\bm{\nu}]+ \mathcal{O}(\omega\mu) \\
    &= -2\ii\omega\sqrt{\frac{\rho}{\lambda+2\mu}}\theta_2(-\frac{1}{2}\bm{I} + (\widetilde{\mathcal{K}}_+^{\alpha,0})^*)^{-1}\left((-\frac{1}{2}\bm{I}+(\mathcal{K}_+^{\alpha,0})^*)^{-1}[\bm{h(\theta)}] - \bm{h(\theta)} - \sum\limits_{i=1}^3 (x_2\bm{\theta}, \bm{\psi}^{(i)})\bm{f}^{(i)}\right)\\&+ 2\lambda\omega\sqrt{\frac{\rho}{\lambda+2\mu}}(-\frac{1}{2}\bm{I}+(\widetilde{\mathcal{K}}_+^{\alpha,0})^*)^{-1}[\bm{\nu}]+\mathcal{O}(\omega\mu).
  \end{align*}
  The third equality holds because
  $$
    (-\frac{1}{2}\bm{I} + (\mathcal{K}_+^{\alpha,0})^*)(\mathcal{S}_+^{\alpha,0})^{-1}[\bm{f}^{(i)}] = 0.
  $$
  For the sake of brevity, we define
  \begin{align*}
    &\mathscr{F}_{\perp,1} := \\&-2\ii\omega\sqrt{\frac{\rho}{\lambda+2\mu}}\theta_2\left((-\frac{1}{2}\bm{I}+(\mathcal{K}_+^{\alpha,0})^*)^{-1}[\bm{h(\theta)}] - \bm{h(\theta)} - \sum\limits_{i=1}^3 (x_2\bm{\theta}, \bm{\psi}^{(i)})\bm{f}^{(i)}\right)+2\lambda\omega\sqrt{\frac{\rho}{\lambda+2\mu}}\bm{\nu},
  \end{align*}
  then
  $$
    (\mathcal{A}_0^{-1}\bm{F}_\perp)_1 = (-\frac{1}{2}\bm{I} + (\widetilde{\mathcal{K}}_+^{\alpha,0})^*)^{-1}[\mathscr{F}_{\perp,1}] + \mathcal{O}(\omega\mu).
  $$
  Similarly, we have
  \begin{align*}
    (\mathcal{A}_0^{-1}\bm{F}_\perp)_2 = -(\mathcal{S}_+^{\alpha,0})^{-1}[\bm{F}_{\perp, 1}] = -2\ii\omega\sqrt{\frac{\rho}{\lambda+2\mu}}\theta_2(\mathcal{S}_+^{\alpha,0})^{-1}[x_2\bm{\theta}] +\mathcal{O}(\omega\mu^2).
  \end{align*}
  \item \emph{$\tau_i\widetilde{\mathcal{A}}_0^{-1}\mathcal{A}_1[\Psi_i]$ term}.

  From \eqref{eq:a1expansion} we know that the leading order term of $\widetilde{\mathcal{A}}_0^{-1}\mathcal{A}_1[\Psi_i]$ is $\mathcal{O}(\omega^2) + \mathcal{O}(\widetilde{\mu}^{-1})$, therefore the contribution of this term in $\phi$ is minor compared to $\tau_i\bm{\psi}^{(i)}$ from $\Phi_0$. Nevertheless, the contribution to $\bm{\psi}$ should not be neglected, since the blowing $\bm{\tau}$ will cancel the decaying $\widetilde{\mu}^{-1}$. We have
  \begin{align*}
    (\widetilde{\mathcal{A}}_0^{-1}\mathcal{A}_1[\Psi_i])_2 = \frac{1}{\widetilde{\mu}}\sum\limits_{i=1}^3\sum\limits_{j=1}^3\tau_is_{ij}(-\mathcal{S}_+^{\alpha,0})^{-1}[\bm{f}^{(j)}] + \mathcal{O}(\omega^2\mu)= -\frac{\mu}{\widetilde{\mu}}\sum\limits_{i=1}^3\tau_i\bm{\psi}^{(i)} + \mathcal{O}(\omega^2\mu).
  \end{align*}

  \item \emph{$\widetilde{\mathcal{A}}_0^{-1}\mathcal{A}_1\widetilde{\mathcal{A}}_0^{-1}[F]$ term.}

  From \eqref{eq:a1expansion} we know that the leading order term of $\widetilde{\mathcal{A}}_0^{-1}\mathcal{A}_1\widetilde{\mathcal{A}}_0^{-1}[F]$ is even higher order terms compared to $\tau_i\widetilde{\mathcal{A}}_0^{-1}\mathcal{A}_1[\Psi_i]$. Hence we neglect the contribution of this term in both $\bm{\phi}$ and $\bm{\psi}$.
\end{itemize}
\subsubsection{Near-field and far-field estimation}

By summarizing all the results we have got in the previous section, we can conclude that
\begin{align*}
    \begin{pmatrix}\bm{\phi} \\ \bm{\psi}\end{pmatrix} &= \Phi = \Phi_\perp + \sum\limits_{i=1}^3\tau_i\Psi_i \\
    &=\begin{pmatrix}(-\frac{1}{2}\bm{I}+(\widetilde{\mathcal{K}}_+^{\alpha,0})^*)^{-1}[\mathscr{F}_{\perp,1}]+\sum\limits_{i=1}^3\tau_i\bm{\psi}^{(i)}+ \mathcal{O}(\omega^2) + \mathcal{O}(\omega\mu) + \mathcal{O}(\frac{1}{\widetilde{\mu}^2}) \\ -2\ii\omega\sqrt{\frac{\rho}{\lambda+2\mu}}\theta_2(\mathcal{S}_+^{\alpha,0})^{-1}[x_2\bm{\theta}]-\frac{\mu}{\widetilde{\mu}}\sum\limits_{i=1}^3\tau_i\bm{\psi}^{(i)} + \mathcal{O}(\omega^2\mu)+\mathcal{O}(\omega\mu^2) \end{pmatrix} ,
\end{align*}
where $\alpha_i$ are defined in \eqref{eq:alphaexpression}.

To conclude, from \eqref{eq:int_repre} we observe that, under the assumption that $\widetilde{\mu}\rightarrow +\infty$, $\mu\rightarrow 0$, the sub-wavelength resonance frequency $\omega\rightarrow 0$, and the field inside scatterer is determined by
\begin{equation}\label{eq:innerfield}
\begin{aligned}
     \widetilde{\mathcal{S}}_+^{\alpha,\omega}[\bm{\phi}] = -2\ii\omega\sqrt{\frac{\rho}{\lambda+2\mu}}\theta_2\bm{u}_{int} + \sum\limits_{i=1}^3\tau_i\widetilde{\mathcal{S}}_+^{\alpha,\omega}[\bm{\psi}^{(i)}] + \mathcal{O}(\frac{\omega\mu}{\widetilde{\mu}}) + \mathcal{O}(\frac{\omega^2}{\widetilde{\mu}}) + \mathcal{O}(\frac{1}{\widetilde{\mu}^3}).
\end{aligned}
\end{equation}
where $\bm{u}_{int} := \widetilde{\mathcal{S}}_+^{\alpha,0}(-\bm{I}/2+(\widetilde{\mathcal{K}}_+^{\alpha,0})^*)^{-1}[\mathscr{F}_{\perp,1}]$ satisfying the following equation system
\begin{equation*}
\left\{
\begin{aligned}
    &\mathcal{L}^{\widetilde{\lambda}, \widetilde{\mu}}\bm{u}_{int} = 0 \quad \mathrm{in}~D \\
    &\frac{\partial \bm{u}_{int}}{\partial\widetilde{\bm{\nu}}}|_- = \mathscr{F}_{\perp, 1} \quad \mathrm{on}~\partial D.
\end{aligned}
\right.
\end{equation*}

We observe that the field inside scatterer blows up due to the blow up of $\tau_i$ when $\omega$ goes to resonance frequency in Theorem \ref{thm:resonantfrequency}.

From expansion of $\bm{G}_+^{\alpha,\omega}$ and definition of $\mathcal{S}_+^{\alpha,\omega}$, we have $\|\mathcal{S}_+^{\alpha,\omega} - \mathcal{S}_+^{\alpha,0}\| = \mathcal{O}(\omega)$, where the norm is taken in space of bounded operators from $(L^2(\partial D))^2$ to $(L^2(\partial D))^2$. The scattered field outside scatterer is determined by
\begin{equation}\label{eq:outerfield}
\begin{aligned}
    \mathcal{S}_+^{\alpha,\omega}[\bm{\psi}] = -2\ii\omega\sqrt{\frac{\rho}{\lambda+2\mu}}\theta_2x_2\bm{\theta} -\frac{\mu}{\widetilde{\mu}}\sum\limits_{i=1}^3\tau_i\mathcal{S}_+^{\alpha,0}[\bm{\psi}^{(i)}]+ \mathcal{O}(\omega\mu\widetilde{\mu}^{-1})+ \mathcal{O}(\omega^2\mu)+\mathcal{O}(\omega\mu^2).
\end{aligned}
\end{equation}

Now we calculate the far-field approximation of $\mathcal{S}_+^{\alpha,\omega}[\bm{\psi}]$. In the limiting case when $\omega\rightarrow 0$, $l_p=l_s=1$. $l_p$ and $l_s$ can never be $0$ since $\alpha<\mathrm{min}(k_p,k_s)$. As a result, $l=0$ is the only one propagating mode when $\omega\rightarrow 0$. In the far-field scenario, $x_2\rightarrow +\infty$, so $x_2-y_2>0$ and $x_2+y_2>0$. The corresponding leading order term of propagating Green's function is
$$
    \bm{G}^\alpha_{+,1,0}(\bm{x},\bm{y}) = -\frac{y_2}{2\pi\rho}\begin{pmatrix}
    \rho & \delta(\alpha_p-\alpha_s) \\
    \delta(\alpha_p-\alpha_s) & \rho\frac{\mu}{\lambda+2\mu}
    \end{pmatrix}e^{\ii \alpha x_1}.
$$
When $x_2\rightarrow +\infty$, the evanescent part in $\mathcal{S}_+^{\alpha,0}[\bm{\psi}^{(i)}]\rightarrow 0$, so when away from metascreen,
\begin{equation*}
\begin{aligned}
    \mathcal{S}_+^{\alpha,\omega}[\bm{\psi}]&\approx \bm{u}_{sc,\infty}^{\alpha,\omega}:= -2\ii\omega\sqrt{\frac{\rho}{\lambda+2\mu}}\theta_2x_2\bm{\theta} \\&+ \frac{1}{\widetilde{\mu}}\sum\limits_{i=1}^3\tau_i\int_{\partial D}\frac{y_2}{2\pi\rho}\begin{pmatrix}
    \rho & \delta(\alpha_p-\alpha_s) \\
    \delta(\alpha_p-\alpha_s) & \rho\frac{\mu}{\lambda+2\mu}
    \end{pmatrix}e^{\ii \alpha x_1}\bm{\psi}^{(i)}(y)\mathrm{d}\sigma(y)+ \mathcal{O}(\omega^2\mu)+\mathcal{O}(\omega\mu^2).
\end{aligned}
\end{equation*}
As discussed in Remark \ref{rmk:order}, when $\mu = \mathcal{O}(\widetilde{\mu}^{-1})$, we have $\bm{\tau} = \mathcal{O}(\widetilde{\mu})$, and the scattered wave varies greatly when $\omega$ goes to resonant frequency. The main results of this subsection is summarized in the following theorem.
\begin{theorem}
    Under the assumption that $\widetilde{\mu}\rightarrow +\infty$ and $\mu\rightarrow 0$, as the frequency $\omega$ of incident wave approaches sub-wavelength resonance frequency, the field inside each scatterer of elastic metascreen blows up in the order of $\mathcal{O}(\mu^{-1})$. When $\mu = \mathcal{O}(\widetilde{\mu}^{-1})$, the blow-up rate is of order $\mathcal{O}(\widetilde{\mu})$. Furthermore, when subwavelength resonance occurs, the far-field approximation of scattering wave differs greatly when $\mu=\mathcal{O}(\widetilde{\mu}^{-1})$ compared to the non-resonating case.
\end{theorem}


\section{Elastic metascreen with $\mu\rightarrow +\infty$}\label{sec3-3}

In this section we would like to show that the elastic metascreen with background $\mu$ goes to $+\infty$ does not have subwavelength resonant mode. Such material refers to a limit of super-solid material which could hardly slide or twist. To this end we need to calculate asymptotic expression of $\mathcal{A}$ when $\mu \rightarrow +\infty$. Here we assume $\alpha = \delta k_s$, where $0\leq \delta \leq 1$ is a constant. We focus on sub-wavelength resonance throughout this paper, so $\omega$ must go to 0. With those constraints on $\omega$ and $\mu$, $\rho$ and $\alpha$ must satisfy $\alpha/\sqrt{\rho} = \delta\omega/\sqrt{\mu}\rightarrow 0$. In this section we will consider two possibilities, $\alpha\rightarrow 0$ with constant $\rho$, and $\rho\rightarrow\infty$ with constant $\alpha$.  The first case will be discussed in detail, and the result of second case is briefly discussed in Remark \ref{rem:anothercase}. The other possibilities could be discussed in the exact same way. $\widetilde{\rho}$ remains constant in this section.

Keep in mind now $\alpha\rightarrow 0$ while $\rho$ is constant. Similarly as Section \ref{sec:asymptoticgreen}, we asymptotically expand $\gpao$ with respect to $\mu$. By direct calculation, we could verify that $\gpao$ could be expanded in the form of
$$
    \gpao(\bm{x},\bm{y}) =: \sum\limits_{p=1}^{+\infty}\frac{1}{\omega^2}\frac{\omega^{2p}}{\mu^p}\bm{G}^{\alpha,\omega}_{+,p} =:  \sum\limits_{p=1}^{+\infty}\frac{1}{\omega^2}\frac{\omega^{2p}}{\mu^p}\sum\limits_{l\in 2\pi\mathbb{Z}}\bm{G}_{+,p,l}^{\alpha,\omega}(\bm{x},\bm{y}).
$$

We use $\bm{G}_+^{\alpha,\omega}$ to define $\mathcal{S}^{\alpha,\omega}_+$ and $(\mathcal{K}^{\alpha,\omega}_+)^{\ast}$. Following a similar decomposition in \eqref{eq:decomposeA}, we decompose $\mathcal{A}$ as
\begin{equation}\label{eq:decomposeA1}
\mathcal{A} = \mathcal{A}_0+\sum\limits_{p=1}^{+\infty}\omega^{2p-2} (\frac{\omega^2}{\widetilde{\mu}^{p+1}}\mathcal{A}_{p,in} +\frac{1}{\mu^{p}} \mathcal{A}_{p,ex}) =: \mathcal{A}_0+\sum\limits_{p=1}^{+\infty}\omega^{2p-2}\mathcal{A}_p =: \mathcal{A}_0 + B,
\end{equation}
where $\mathcal{A}_{p,in}$, $\mathcal{A}_{p,ex}$ is the expansion of left and right column of $\mathcal{A}$, correspondingly, in the following expressions
\begin{equation}\label{eq:decomposedetail}
\begin{aligned}
    \mathcal{A}_0 &= \begin{pmatrix}
    \widetilde{\mathcal{S}}^{\alpha,0}_+ & 0 \\ -\frac{1}{2}\bm{I} + (\widetilde{\mathcal{K}}^{\alpha,0}_+)^* & -(\frac{1}{2}\bm{I}+(\mathcal{K}^{\alpha,0}_+)^*)
    \end{pmatrix}, \\
    \mathcal{A}_{p,in} &= \begin{pmatrix}
    \widetilde{\mathcal{S}}^\alpha_{+,p+1} & 0 \\ (\widetilde{\mathcal{K}}^\alpha_{+,p+1})^* & 0
    \end{pmatrix}, \qquad
    \mathcal{A}_{p,ex} = \begin{pmatrix}
    0 & -\mathcal{S}^{\alpha}_{+,p} \\ 0 & -\frac{\omega^2}{\widetilde{\mu}}(\mathcal{K}^\alpha_{+,p+1})^*
    \end{pmatrix}, \quad \mathrm{for~} p>0.
\end{aligned}
\end{equation}

Now we calculate the characteristic value of $\mathcal{A}$. Following the steps mentioned at the beginning of Section \ref{sec:resonantfrequency}, we regard $\mathcal{A}_0$ as a perturbation of $\mathcal{A}$ when $\mu\rightarrow \infty$, and study the invertibility of $\mathcal{A}_0$. We assume $(\bm{f}_1, \bm{f}_2)^\top$ is a nontrivial vector satisfying $\mathcal{A}_0[(\bm{f}_1,\bm{f}_2)^\top]=0$, then it immediately follows that
\begin{equation}\label{eq:eqmuinfi}
\left\{
\begin{aligned}
    &\widetilde{\mathcal{S}}_+^{\alpha,0}[\bm{f}_1] = 0, \\
    &(-\frac{1}{2}\bm{I} + (\widetilde{\mathcal{K}}^{\alpha,0}_+)^*)[\bm{f}_1] - (\frac{1}{2}\bm{I}+(\mathcal{K}^{\alpha,0}_+)^*)[\bm{f}_2] = 0.
\end{aligned}
\right.
\end{equation}
Here we still impose Assumption \ref{asump:invertibility} to ensure the invertibility of $\mathcal{S}_+^{\alpha,0}$. Therefore $\bm{f}_1$ must be 0. We now turn to discuss the invertibility of $\frac{1}{2}\bm{I}+(\mathcal{K}^{\alpha,0}_+)^*$. The following lemma proves the uniqueness of solution of quasi-periodic homogeneous Lam\'{e} equation with Neumann boundary on $\partial D$ and Dirichlet boundary on $\mathbb{R}_+^2$.

\begin{lemma}\label{lem:quasiunique}
Let $S_L = \{\bm{x}=(x_1,x_2): \,0\leq x_1\leq L,\; x_2\in\mathbb{R}\}$, $\bm{u}=\bm{u}(\bm{x})$ be defined in $S_L\backslash D$ satisfying
$$
\mathcal{L}^{\lambda,\mu} \bm{u} = 0 \quad \mathrm{in }~~ S_L\setminus \overline{D},
$$
and
$$
\frac{\partial \bm{u}}{\partial \bm{\nu}} = 0 \quad \mathrm{on } ~~\partial D,\qquad \bm{u} = 0 \quad \mathrm{for~~ } x_2=0.
$$
Furthermore, assume that $\bm{u}$ is quasi-periodic in the $x_1$-direction, i.e.
$$
\bm{u}(x_1+L,x_2) = e^{-i\alpha_1 L}\bm{u}(x_1,x_2).
$$
Then $\bm{u}\equiv 0$.
\end{lemma}
\begin{proof}

Define
$$
S_L' = S_L \setminus \overline{D}.
$$
Multiplying the equation
$$
\mathcal{L}^{\lambda,\mu} \bm{u} = 0 \quad \text{in } S_L'
$$
by the complex conjugate $\overline{\bm{u}}$ and integrating over $S_L'$ yields, after an integration by parts,
$$
 \int_{S_L'} \lambda|\nabla\cdot \bm{u}|^2 + 2\mu \nabla^s \bm{u}:\nabla^s\overline{\bm{u}}\mathrm{d}x = \int_{\partial S_L'} \overline{\bm{u}}\cdot \frac{\partial \bm{u}}{\partial\bm{\nu}}\mathrm{d}s = (\int_{\partial D} + \int_{S_0} + \int_{S_L}+\int_{\mathbb{R}_{+,L}}) \overline{\bm{u}}\cdot \frac{\partial \bm{u}}{\partial\bm{\nu}}\mathrm{d}s,
$$
where
\begin{align*}
S_0 := \{\bm{x}\in  S_L: x_1=0\}, && S_L := \{\bm{x}\in S_L: x_1=L\},&& \mathbb{R}_{+,L} := \{\bm{x}\in S_L: x_2=0\}.
\end{align*}
denotes the left, right cell boundary and $x_1$ axis, correspondingly. On $\partial D$, the Neumann condition implies zero contribution. On $\mathbb{R}_{+,L}$, the Dirichlet condition also implies zero contribution. On the left and right cell boundary, the outward normal vector on $S_0$ and $S_L$ differs by a negative sign. By quasi-periodicity, $\bm{u}(L,x_2) = e^{-\ii\alpha L}\bm{u}(0,x_2)$, and the derivatives of $\bm{u}$ differs by a factor of $e^{-\ii\alpha L}$ on $S_0$ and $S_L$ as well. Therefore, when calculating $\displaystyle \overline{\bm{u}}\cdot\frac{\partial \bm{u}}{\partial\bm{\nu}}$ on $S_L$, the factor $e^{-\ii\alpha L}$ brought in by $\displaystyle \frac{\partial \bm{u}}{\partial\bm{\nu}}$ gets cancelled with the factor brought in by $\overline{\bm{u}}$, which means
$$
    \int_{S_L} \overline{\bm{u}(L,x')}\cdot \frac{\partial \bm{u}}{\partial \bm{\nu}}(L,x')\mathrm{d}s = -\int_{S_0} \overline{\bm{u}(L,x')}\cdot \frac{\partial \bm{u}}{\partial \bm{\nu}}(L,x')\mathrm{d}s.
$$
Therefore,
$$
    \int_{S_L'} \lambda|\nabla\cdot \bm{u}|^2 + 2\mu \nabla^s \bm{u}:\nabla^s\overline{\bm{u}}\mathrm{d}x = 0.
$$
From Section 1.4.2 of \cite{HabibElastic} we know that the quadratic form
$$
Q(\bm{u},\bm{u}) := \int_{S_L'} \lambda|\nabla\cdot \bm{u}|^2 + 2\mu \nabla^s \bm{u}:\nabla^s\overline{\bm{u}}\mathrm{d}x
$$
is positive definite. This implies $\bm{u}$ must be a rigid body motion, i.e., $\bm{u}(\bm{x}) = \bm{Ax + b}$, where $\bm{A}$ is a skew-symmetric matrix and $\bm{b}$ is a constant vector. Nevertheless, the quasi-periodicity condition
$$
\bm{u}(x_1+L,x_2) = e^{-i\alpha_1 L}\,\bm{u}(x_1,x_2)
$$
imposes a nontrivial phase shift that is incompatible with any nonzero rigid body motion. Thus, the only possibility is that $\bm{u} \equiv 0$ in $S_L'$.
\end{proof}
\begin{proposition}\label{prop:non_scatter}
When $\alpha\neq 0$, $\frac{1}{2}\bm{I}+(\mathcal{K}_+^{\alpha,0})^{\ast}:L^2(\partial D)^2\to L^2(\partial D)^2$ is invertible.
\end{proposition}
\begin{proof}
We first show the injectivity of $\frac{1}{2}\bm{I}+(\mathcal{K}_+^{\alpha,0})^{\ast}$.
Suppose $\bm{\phi}\in L^2(\partial D)^2$ satisfies $\big(\frac{1}{2}\bm{I}+(\mathcal{K}_+^{\alpha,0})^{\ast}\big)\bm{\phi}=0$ on $\partial D$. Using notations defined in Lemma \ref{lem:quasiunique}, $\bm{u}:=\mathcal{S}_+^{\alpha,0}[\bm{\phi}]$ satisfies
\begin{equation}\label{eq:lemma1}
\left\{
	\begin{array}{ll}
	\mathcal{L}^{\lambda, \mu}\bm{u}=\bm{0}	,  &   \text{in} \ S_L\setminus\bar{D}, \medskip \\
                 \displaystyle  \left.\frac{\partial \bm{u}}{\partial \bm{\nu}}\right\vert_+ =\bm{0},   &\text{on} \   \partial D ,\medskip \\
                 \bm{u} = 0 , & \text{for} \ x_2=0, \medskip \\
                 \bm{u} \ \text{is}\ \alpha\text{-quasi-periodic}.
	\end{array}
	\right.
\end{equation}
From Lemma \ref{lem:quasiunique}, we have the uniqueness of the solution of \eqref{eq:lemma1}. Therefore $\bm{u}=0$ is the unique solution and by the uniqueness theorem inside $D$ (See for example \cite{DahlbergKenigVerchota1988}), $\bm{u}$ is also constant in $D$. Therefore, $\displaystyle \bm{\phi}=\left.\frac{\partial \bm{u}}{\partial \bm{\nu}}\right\vert_+-\left.\frac{\partial \bm{u}}{\partial \bm{\nu}}\right\vert_- =\bm{0} $.
Notice that compactness of $(\mathcal{K}_+^{\alpha,0})^*$ is a standard result, which can be proved in a similar manner like in \cite{ColtonKress2013}. The invertibility of $\frac{1}{2}\bm{I}+(\mathcal{K}_+^{\alpha,0})^{\ast}$ then follows directly from Fredholm alternative.
\end{proof}

With our assumption on $D$ and Proposition \ref{prop:non_scatter}, we have the following conclusion.

\begin{theorem}\label{thm:dirichletresonance}
    Let $D$ be chosen such that $\mathcal{S}_+^{\alpha,0}$ is invertible, then $\mathcal{A}_0$ is invertible. That is, when $\mu\rightarrow \infty$, there is no sub-wavelength resonance in system \eqref{eq:elas_sys}.
\end{theorem}

\begin{remark}
It is worth mentioning here that there are still resonances in system \eqref{eq:elas_sys}. It is straightforward to observe that when $\mu\rightarrow\infty$ and $\omega$ is not very close to 0, we could expand $\mathcal{A}$ with respect to $\mu$ and get leading order term
$$
    \mathcal{A}_0 = \begin{pmatrix}
    \widetilde{\mathcal{S}}^{\alpha,\omega}_+ & 0 \\ -\frac{1}{2}\bm{I} + (\widetilde{\mathcal{K}}^{\alpha,\omega}_+)^* & -(\frac{1}{2}\bm{I}+(\mathcal{K}^{\alpha,0}_+)^*)
    \end{pmatrix}.
$$
In this case, the invertibility of $\mathcal{A}_0$ solely depends on the invertibility of $\widetilde{\mathcal{S}}^{\alpha,\omega}_+$, which fails to be invertible when $\omega$ matches the eigenvalue of interior Dirichlet problem. This resonance, however, is not sub-wavelength resonance. The details of discussion of this resonance type could be found in \cite{HabibARMA}.
\end{remark}

\begin{remark}\label{rem:anothercase}
In this remark we briefly discuss the case when $\alpha$ is constant with $\rho\rightarrow\infty$, which also satisfies $\alpha/\sqrt{\rho}\rightarrow 0$. In this case we have $\rho\omega^2/\mu = \alpha^2/\delta^2$, therefore
$$
    \gamma_{l,s} = \sqrt{|(\alpha+l)^2-\frac{\alpha^2}{\delta^2}|}, \quad \gamma_{l,p} = \sqrt{|(\alpha+l)^2-\frac{\alpha^2}{\delta^2}C_{\lambda,\mu}|}.
$$
This implies that $\bm{G}^{\alpha,\omega}$ is of order $O(1/\mu)$, and $\mathcal{S}_+^{\alpha,\omega}$ has no $O(1)$ term. We then get exactly the same $\mathcal{A}_0$ as in \eqref{eq:decomposedetail}. The same resonance condition as Theorem \ref{thm:dirichletresonance} follows.
\end{remark}


\section{Conclusion}\label{sec13}

In this paper we studied the criterion of subwavelength resonance of one-dimensional elastic metascreen in two-dimensional space. Using the quasi-periodic Green's functions of Lam\'{e} system and corresponding layer potential operators calculated in Section \ref{sec3}, we proved that the subwavelength resonance will emerge when the inclusion shear modulus $\widetilde{\mu}$ tends to infinity, and the quotient of background shear modulus $\mu$ and inclusion density $\widetilde{\rho}$ goes to zero, simultaneously. This corresponds to the high contrast scenario $\widetilde{\mu}/\mu\rightarrow \infty$ discussed in other references. We further proved that no subwavelength resonance will appear when the background shear modulus $\mu$ tends to infinity, no matter how the other parameters change. This corresponds to the high contrast scenario $\widetilde{\mu}/\mu\rightarrow 0$. The reason that why we only discussed the high contrast between $\widetilde{\mu}$ and $\mu$ is, that as far as we are concerned, the literatures on devising elastic metascreen and metamaterial all required high contrast in shear modulus plus some additional requirements. It is worth mentioning that with exactly the same methods, we could discuss whether or not the subwavelength resonance exists in the $(\lambda, \widetilde{\lambda})$ high contrast pair or the $(\rho, \widetilde{\rho})$ high contrast pair, and calculate the subwavelength resonance frequency in case it exists. The quasi-periodic Green's function and occurrence condition of subwavelength resonance for the one-dimensional and two-dimensional elastic metascreen in three-dimensional space could be achieved following this paper without much difficulty.

\backmatter





\bmhead{Acknowledgements}

W. W. was supported by NSFC grant (12301539), NSFC/RGC Joint Research Fund (project N\_CityU101/21) and youth growth grant of Department of Science and Technology of Jilin Province (20240602113RC). Y. H. was supported by Shenzhen Science and Technology Program (Grant No. RCBS20231211090750090) and NSFC grant (12401564).

\bmhead{Data availability} No datasets were generated or analysed during the current study, because our work proceeds within a theoretical
 and mathematical approach.

\section*{Declarations}


\bmhead{Conflict of interest} The author has no conflict of interest to declare that are relevant to the content of this article.








\begin{appendices}

\section{Characteristic value}\label{sec:appa}

For the sake of completeness, we introduce Rouch\'{e}'s theorem and all related concepts here. This part is an excerpt from Chapter 1 of \cite{HabibLayer}.

Let $\mathcal{B}, \mathcal{B}'$ be Banach spaces. Let $\mathfrak{U}(z_0)$ be the set of all operator-valued functions with values in $\mathcal{L}(\mathcal{B}, \mathcal{B}')$ which are holomorphic in some neighborhood of $z_0$, except possibly at $z_0$. The point $z_0$ is called a \emph{characteristic value} of $A(z) \in \mathfrak{U}(z_0)$ if there exists a vector-valued function $\phi(z)$ with values in $\mathcal{B}'$ such that
\begin{itemize}
\item  $\phi(z)$ is holomorphic at $z_0$ and $\phi(z_0) = 0$,
\item  $A(z)\phi(z)$ is holomorphic at $z_0$ and vanishes at this point.
\end{itemize}
Here, $\phi(z)$ is called a \emph{root function} of $A(z)$ associated with the characteristic value $z_0$. The vector $\phi_0 = \phi(z_0)$ is called an eigenvector. The closure of the linear set of eigenvectors corresponding to $z_0$ is denoted by $\mathrm{Ker}A(z_0)$.

Suppose that $z_0$ is a characteristic value of the function $A(z)$ and $\phi(z)$ is an associated root function. Then there exists a number $m(\phi) \geq 1$ and a vector-valued function $\psi(z)$ with values in $\mathcal{B}'$, holomorphic at $z_0$, such that
$$
A(z)\phi(z) = (z - z_0)^{m(\phi)}\psi(z),\quad \psi(z_0) \neq 0.
$$
The number $m(\phi)$ is called the \emph{multiplicity} of the root function $\phi(z)$.

For $\phi_0 \in \mathrm{Ker}A(z_0)$, we define the rank of $\phi_0$, denoted by $\mathrm{rank}(\phi_0)$, to be the maximum of the multiplicities of all root functions $\phi(z)$ with $\phi(z_0) = \phi_0$. Suppose that $n = \mathrm{dim} \mathrm{Ker}A(z_0) < +\infty$ and that the ranks of all vectors in $\mathrm{Ker}A(z_0)$ are finite. A system of eigenvectors $\phi_0^j, j=1,\cdots, n$ is called a canonical system of eigenvectors of $A(z)$ associated to $z_0$ if their ranks possess the following property: for $j = 1,\cdots, n$, $\mathrm{rank}(\phi^j_0)$ is the maximum of the ranks of all eigenvectors in the direct complement in $\mathrm{Ker}A(z_0)$ of the linear span of the vectors $\phi_0^1, \cdots, \phi_0^{j-1}$. We call
$$
N(A(z_0)) :=\sum\limits_{j=1}^n \mathrm{rank}(\phi_0^j)
$$
the \emph{null multiplicity} of the characteristic value $z_0$ of A(z). If $z_0$ is not a characteristic value of $A(z)$, we put $N(A(z_0)) = 0$.

Suppose that $A^{-1}(z)$ exists and is holomorphic in some neighborhood of $z_0$, except possibly at $z_0$. Then the number
$$
M(A(z_0)) = N(A(z_0)) - N(A^{-1}(z_0))
$$
is called the \emph{multiplicity} of $z_0$. If $z_0$ is a characteristic value and not a pole of $A(z)$, then $M(A(z_0)) = N(A(z_0))$ while $M(A(z_0)) = -N(A^{-1}(z_0))$ if $z_0$ is a pole and not a characteristic value of A(z).

Suppose that $z_0$ is a pole of the operator-valued function $A(z)$ and the Laurent series expansion of $A(z)$ at $z_0$ is given by
\begin{equation}\label{eq:normalpoints}
A(z) = \sum\limits_{j\geq -s}(z - z_0)^jA_j.
\end{equation}
If the operators $A_{-j}, j = 1, \cdots, s$, have finite-dimensional ranges, then $A(z)$ is called \emph{finitely meromorphic} at $z_0$.The operator-valued function $A(z)$ is said to be of Fredholm type (of index zero) at the point $z_0$ if the operator $A_0$ in \eqref{eq:normalpoints} is Fredholm (of index zero). If $A(z)$ is holomorphic and invertible at $z_0$, then $z_0$ is called a \emph{regular point} of $A(z)$. The point $z_0$ is called a \emph{normal point} of $A(z)$ if $A(z)$ is finitely meromorphic, of Fredholm type at $z_0$, and regular in a neighborhood of $z_0$ except at $z_0$ itself.

Let $V$ be a simply connected bounded domain with rectifiable boundary $\partial V$. An operator-valued function $A(z)$ which is finitely meromorphic and of Fredholm type in $V$ and continuous on $\partial V$ is called \emph{normal} with respect to $\partial V$ if the operator $A(z)$ is invertible in $V$, except for a finite number of points of $V$ which are normal points of $A(z)$.

\begin{theorem}[Generalized Rouch\'{e}'s theorem]\label{thm:rouche}
Let $A(z)$ be an operator-valued function which is normal with respect to $\partial V$. If an operator-valued function $S(z)$ which is finitely meromorphic in $V$ and continuous on $\partial V$ satisfies the condition
$$
    \|A^{-1}(z)S(z)\|_{\mathcal{L}(\mathcal{B},\mathcal{B})}<1, \quad z\in\partial V,
$$
then $A(z) + S(z)$ is also normal with respect to $\partial V$ and
$$
\mathcal{M}(A(z); \partial V) = \mathcal{M}(A(z) + S(z); \partial V ).
$$
\end{theorem}

\end{appendices}


\bibliography{sn-article}

\end{document}